# SINGULAR PRINCIPAL BUNDLES OVER HIGHER DIMENSIONAL MANIFOLDS AND THEIR MODULI SPACES[†]

### ALEXANDER H.W. SCHMITT

ABSTRACT. In this note, we introduce the notion of a singular principal $G$-bundle, associated to a reductive algebraic group $G$ over the complex numbers by means of a faithful representation $\rho' \colon G \longrightarrow \mathrm{SL}(V)$. This concept is meant to provide an analogon to the notion of a torsion free sheaf as a generalization of the notion of a vector bundle.

We will construct moduli spaces for these singular principal bundles which compactify the moduli spaces of stable principal bundles.

## INTRODUCTION

Given an algebraic group $G$, one has the notion of a principal $G$-bundle over an algebraic variety $X$. If $G = \mathrm{GL}(r)$, this notion is equivalent to the notion of a vector bundle of rank $r$. In the past, there has been enormous research in the theory of vector bundles. Large part of the theory focusses on semistable vector bundles. A major achievement is the construction of the moduli space of semistable vector bundles over curves, due to Seshadri. Its generalizations, obtained by Gieseker and Maruyama, reveal that, on higher dimensional manifolds, one has to include semistable torsion free sheaves as well in order to obtain projective moduli spaces. The same goes for singular curves. Although the theory of principal $G$-bundles for other reductive groups $G$ has received growing interest from mathematical physics (see [3] and [2]), the above findings were only partially generalized. Ramanathan gave in [13] an intricate GIT construction of the moduli space of semistable $G$-bundles on a smooth projective curve, proceeding in several steps, according to the chain of homomorphisms

$$G \longrightarrow G/\mathrm{Cen}(G) \overset{\mathrm{Ad}}{\cong} \mathrm{Aut}^0(\mathrm{Lie}(G)) \hookrightarrow \mathbb{C}^* \times \mathrm{Aut}(\mathrm{Lie}(G)) \hookrightarrow \mathrm{GL}(\mathrm{Lie}(G)).$$

Hyeon [8] has generalized Ramanathan's construction to give the moduli spaces of stable principal bundles over higher dimensional base schemes, but the resulting moduli spaces are only quasi-projective, and the constructions do not suggest any natural compactification. The necessary singular generalizations of principal $G$-bundles have been considered only in the case of the classical groups $\mathrm{O}(r, \mathbb{C})$, $\mathrm{SO}(r, \mathbb{C})$, and $\mathrm{Sp}(r, \mathbb{C})$. For these groups, the principal $G$-bundles have natural interpretations as vector bundles with additional structures, and these can be easily extended to the setting of torsion free sheaves [4].

---

In this note, we will propose a more general approach, working for arbitrary reductive groups. Given a faithful representation $\rho' \colon G \longrightarrow \mathrm{SL}(V)$, we will introduce the notion of a singular principal bundle which mimics the construction of a principal $G$-bundle from a principal $\mathrm{GL}(V)$-bundle $\mathcal{F}$ and a section $X \longrightarrow \mathcal{F}/G$. This concept can again be formulated entirely in the realm of torsion free sheaves with an additional structure which, however, looks more complicated than usual.

Luckily, for the GIT construction, these objects are not much harder to treat than the tensor fields considered in [15] and [4]. Therefore, we obtain projective moduli spaces for singular $G$-principal bundles. Our construction is more direct than the one used by Ramanathan and Hyeon.

After finishing a former version of this paper, Gómez and Sols [5] announced the construction of projective moduli spaces for so called principal $G$-sheaves where $G$ is a reductive group. We will discuss the relation with our approach at the end of this paper.

**Notation.** $X$ will be a projective manifold over the field of complex numbers, $\mathcal{O}_X(1)$ will be a fixed ample line bundle on $X$. This defines, for every coherent $\mathcal{O}_X$-module $\mathcal{A}$, the Hilbert polynomial $P(\mathcal{A})$ with $P(\mathcal{A}(n)) = \chi(\mathcal{A}(n))$ for all $n \geq 0$ as well as the slope $\mu(\mathcal{A}) := \bigl(c_1(\mathcal{A}).c_1(\mathcal{O}_X(1))^{\dim X - 1}[X]\bigr)/\operatorname{rk}\mathcal{A}$ and $\mu_{\max}(\mathcal{A}) := \max\{\, \mu(\mathcal{A}') \mid 0 \subsetneq \mathcal{A}' \subset \mathcal{A} \,\}$.

**Acknowledgment.** I am indebted to the anonymous and patient referees for pointing out several major and minor inaccuracies. Their help was crucial for bringing the paper into its final form.

## 1. Preliminaries

1.1. **Reminder on Spec.** Let $X$ be a scheme and $\mathcal{B}$ a quasi-coherent sheaf of $\mathcal{O}_X$-algebras. This provides the affine morphism $\pi \colon \underline{\operatorname{Spec}}(\mathcal{B}) \longrightarrow X$ which has the following base change property: For any morphism $f \colon X' \longrightarrow X$, the morphism $\pi'$ in the cartesian diagram

$$
\begin{array}{ccc}
X' \times_X \underline{\operatorname{Spec}}(\mathcal{B}) & \longrightarrow & \underline{\operatorname{Spec}}(\mathcal{B}) \\
\pi' \downarrow & & \downarrow \pi \\
X' & \xrightarrow{\ \ f\ \ } & X
\end{array}
$$

is again affine, and $X' \times_X \underline{\operatorname{Spec}}(\mathcal{B})$ is canonically isomorphic to $\underline{\operatorname{Spec}}(f^*\mathcal{B})$ via a canonical isomorphism $\pi'_* \mathcal{O}_{X' \times_X \underline{\operatorname{Spec}}(\mathcal{B})} \longrightarrow f^*\mathcal{B}$ of $\mathcal{O}_{X'}$-algebras, [6], (1.5.2), p. 12.

1.2. **Universal good quotients.** Let $G$ be an algebraic group acting on the scheme $X$. A *universal good quotient* is a pair $(Y, \varphi)$, consisting of a scheme $Y$ and a $G$-invariant morphism $\varphi \colon X \longrightarrow Y$, such that for every morphism $Y' \longrightarrow Y$, the induced morphism $\varphi' \colon Y' \times_Y X \longrightarrow Y'$ is a good quotient for the induced $G$-action. Recall that Mumford's GIT [11] provides universal good quotients.

1.3. **Reynolds operators.** Let $\rho \colon G \longrightarrow \mathrm{GL}(V)$ be a locally finite representation of the reductive group $G$ on the possibly infinite dimensional vector space $V$. This means that every element $v \in V$ is contained in a finite dimensional $G$-invariant subspace of $V$. Let $V^G \subset V$ be the subspace of those $v$ for which $\rho(g)(v) = v$ for all $g \in G$. Then, there is a unique map $\mathcal{R} \colon V \longrightarrow V$ of $G$-modules with $\mathcal{R}^2 = \mathcal{R}$ and $\mathcal{R}(v) = v$ if and only if $v \in V^G$. The map $\mathcal{R}$ is called the *Reynolds operator*. It commutes with maps between $G$-modules (see [11], p. 26f).

**1.4. A functoriality property.** We can use the results of Section 1.1 - 1.3 to deal with the following situation: Let $G$ be a reductive algebraic group, acting trivially on the quasi-projective scheme $X$, and let $\mathcal{A}$ be a $G$-linearized coherent $\mathcal{O}_X$-module. Then, the algebra $\mathcal{B} := \mathrm{Sym}^* \mathcal{A}$ inherits a $G$-linearization, providing a $G$-action on $\underline{\mathrm{Spec}}(\mathcal{B})$ ([7], p. 54). We can form the GIT-quotient $\underline{\mathrm{Spec}}(\mathcal{B}) /\!\!/ G$. For this, one chooses an ample line bundle $\mathcal{L}$ on $X$ and works with the trivial linearization in $\pi^*\mathcal{L}$. One verifies

$$\underline{\mathrm{Spec}}(\mathcal{B}) /\!\!/ G \quad = \quad \underline{\mathrm{Spec}}(\mathcal{B}^G),$$

where $\mathcal{B}^G$ is the sheaf of $G$-invariant sections of $\mathcal{B}$. By Section 1.1, for any morphism $f\colon X' \longrightarrow X$ the squares in the diagram

$$
\begin{array}{ccc}
\underline{\mathrm{Spec}}(f^*\mathcal{B}) & \longrightarrow & \underline{\mathrm{Spec}}(\mathcal{B}) \\
\downarrow & & \downarrow \\
\underline{\mathrm{Spec}}(f^*(\mathcal{B}^G)) & \longrightarrow & \underline{\mathrm{Spec}}(\mathcal{B}^G) \\
\downarrow & & \downarrow \\
X' & \overset{f}{\longrightarrow} & X
\end{array}
$$

are cartesian. From Section 1.2, we infer

$$f^*(\mathcal{B}^G) \quad = \quad (f^*\mathcal{B})^G.$$

For every affine open set $U = \mathrm{Spec}\, R \subset X$, $\mathcal{B}(U)$ is a finitely generated $\mathbb{C}$-algebra, and the representation of $G$ on $\mathcal{B}(U)$ is locally finite. Since $R \subset \mathcal{B}(U)^G$, the Reynolds operator $\mathcal{R}_U\colon \mathcal{B}(U) \longrightarrow \mathcal{B}(U)^G$ is an $R$-module homomorphism. Note that, by uniqueness of the Reynolds operator, the $\mathcal{R}_U$ glue to a homomorphism of $\mathcal{O}_X$-modules(!) $\mathcal{R}\colon \mathcal{B} \longrightarrow \mathcal{B}^G$. This Reynolds operator commutes with base change, i.e., $\mathcal{R}\colon f^*\mathcal{B} \longrightarrow (f^*\mathcal{B})^G$ coincides with $f^*\mathcal{R}\colon f^*\mathcal{B} \longrightarrow f^*(\mathcal{B}^G)$.

**1.5. A useful construction.** The following result is due to Gómez and Sols ([4], Lemma 3.1).

**Proposition 1.5.1.** *Let $S$ be a noetherian scheme, $\mathcal{A}^1_S$ and $\mathcal{A}^2_S$ coherent sheaves on $S \times X$, and $\varphi_S\colon \mathcal{A}^1_S \longrightarrow \mathcal{A}^2_S$ a homomorphism. Assume that $\mathcal{A}^2_S$ is $S$-flat. Then, there is a closed subscheme $\mathfrak{Y} \subset S$ the closed points of which are those $s \in S$ for which $\varphi_{S|\{s\}\times X} \equiv 0$. More precisely, it has the property that any morphism $f\colon T \longrightarrow S$ factors through $\mathfrak{Y}$, if and only if $(f \times \mathrm{id}_X)^*\varphi_S$ is the zero homomorphism.*

*Proof.* We may write $\mathcal{A}^1_S$ as the quotient of a vector bundle $\mathcal{V}_S$. Let $\varphi'_S\colon \mathcal{V}_S \longrightarrow \mathcal{A}^1_S \longrightarrow \mathcal{A}^2_S$ be the composed homomorphism. This is a homomorphism between $S$-flat sheaves. Choose $n$ large enough, so that the higher direct image sheaves of $\mathcal{V}_S \otimes \pi_X^*\mathcal{O}_X(n)$ and $\mathcal{A}^2_S \otimes \pi_X^*\mathcal{O}_X(n)$ under $\pi_S$ vanish and the evaluation maps $\pi_S^*\pi_{S*}\bigl(\mathcal{V}_S \otimes \pi_X^*\mathcal{O}_X(n)\bigr) \longrightarrow \mathcal{V}_S \otimes \pi_X^*\mathcal{O}_X(n)$ and $\pi_S^*\pi_{S*}\bigl(\mathcal{A}^2_S \otimes \pi_X^*\mathcal{O}_X(n)\bigr) \longrightarrow \mathcal{A}^2_S \otimes \pi_X^*\mathcal{O}_X(n)$ become surjective. Then, $\pi_{S*}\bigl(\varphi'_S \otimes \mathrm{id}_{\pi_X^*\mathcal{O}_X(n)}\bigr)$ is a map between vector bundles, and we may define $\mathfrak{Y}$ as the zero scheme of this map.

Given a morphism $f\colon T \longrightarrow S$, it will factor through the scheme $\mathfrak{Y}$, if and only if $f^*\bigl(\pi_{S*}(\varphi'_S \otimes \mathrm{id}_{\pi_X^*\mathcal{O}_X(n)})\bigr)$ is identically zero. By base change, this is equivalent to the vanishing of $\pi_{T*}\big((f \times \mathrm{id}_X)^*(\varphi'_S \otimes \mathrm{id}_{\pi^*_X \mathcal{O}_X(n)})\big)$. This is in turn equivalent to the vanishing of $\varphi'_T := (f \times \mathrm{id}_X)^*\big(\varphi'_S \otimes \mathrm{id}_{\pi^*_X \mathcal{O}_X(n)}\big)$, by our choice of $n$. The latter homomorphism factorizes as

$$(f \times \mathrm{id}_X)^* \mathcal{V}_S \longrightarrow (f \times \mathrm{id}_X)^* \mathcal{A}^1_S \xrightarrow{\varphi_T} (f \times \mathrm{id}_X)^* \mathcal{A}^2_S$$

with $\varphi_T := (f \times \mathrm{id}_X)^*\big(\varphi_S \otimes \mathrm{id}_{\pi^*_X \mathcal{O}_X(n)}\big)$. As $(f \times \mathrm{id}_X)^* \mathcal{V}_S \longrightarrow (f \times \mathrm{id}_X)^* \mathcal{A}^1_S$ is still surjective, the vanishing of $\varphi'_T$ is equivalent to the one of $\varphi_T$, and we are done. $\square$

**1.6. An extension property.** Let $X$ be a smooth projective manifold, $S$ a noetherian scheme, and $\mathcal{A}_S$ an $S$-flat family of coherent $\mathcal{O}_X$-modules. Let $\mathcal{U}_{\mathrm{lf}}$ be the maximal open subset over which $\mathcal{A}_S$ is locally free. Then, for any $s \in S$, the set $\mathcal{U}_{\mathrm{lf}} \cap \{s\} \times X$ is the maximal open subset where $\mathcal{A}_{S|\{s\} \times X}$ is locally free ([7], Lem. 2.1.7, p. 35).

**Proposition 1.6.1.** *In the above situation, assume* $\mathrm{codim}\big(X \setminus (\mathcal{U}_{\mathrm{lf}} \cap \{s\} \times X), X\big) \geq 2$ *for all $s \in S$. Then, any homomorphism of $\mathcal{O}_{\mathcal{U}_{\mathrm{lf}}}$-modules*

$$\iota \colon \mathcal{A}_{S|U} \longrightarrow \mathcal{O}_{\mathcal{U}_{\mathrm{lf}}}$$

*can be extended in a unique way to a homomorphism of $\mathcal{O}_X$-modules*

$$\iota_* \colon \mathcal{A}_S \longrightarrow \mathcal{O}_X.$$

*Proof.* The arguments used by Maruyama ([10], p. 111f) show that, under the above circumstances, $\mathcal{E}_S \longrightarrow j_*(\mathcal{E}_{S|\mathcal{U}_{\mathrm{lf}}})$ is an isomorphism for any locally free sheaf $\mathcal{E}_S$ on $S \times X$, $j \colon \mathcal{U}_{\mathrm{lf}} \longrightarrow S \times X$ being the inclusion. Thus, the desired extension can be obtained as

$$\mathcal{A}_S \longrightarrow j_*(\mathcal{A}_{S|\mathcal{U}_{\mathrm{lf}}}) \xrightarrow{j_*(\iota)} j_*(\mathcal{O}_{\mathcal{U}_{\mathrm{lf}}}) = \mathcal{O}_X.$$

To see that the extension is unique, note that we can write $\mathcal{A}_S$ as the quotient of a vector bundle $\mathbb{E}_S$, and that $j_*(\mathrm{Hom}(\mathbb{E}_{S|\mathcal{U}_{\mathrm{lf}}}, \mathcal{O}_{\mathcal{U}_{\mathrm{lf}}})) = j_*(\mathbb{E}^\vee_{S|\mathcal{U}_{\mathrm{lf}}}) = \mathbb{E}^\vee_S = \mathrm{Hom}(\mathbb{E}_S, \mathcal{O}_X)$, by the same token. Thus, the composed homomorphism $\mathbb{E}_S \longrightarrow \mathcal{A}_S \longrightarrow \mathcal{O}_S$ is uniquely determined, whence so is $\iota_*$. $\square$

## 2. The moduli problem

**2.1. Principal bundles and their generalizations.** Let $G$ be a reductive algebraic group and $\rho' \colon G \longrightarrow \mathrm{SL}(V)$ a faithful representation.

*Example* 2.1.1. i) If the center of $G$ is trivial, we may take the adjoint representation.

ii) Since semisimple groups are perfect, i.e., $G = (G, G)$ ([1], 14.2, Corollary, p. 182), they do not have any non-trivial characters. Hence, any representation of a semisimple group automatically lands in the special linear group.

iii) If $\rho \colon G \longrightarrow \mathrm{GL}(V)$ is any representation of a reductive group $G$, we get the character $\chi := \det \circ \rho$. Then, $\rho \oplus \chi^{-1} \colon G \longrightarrow \mathrm{GL}(V \oplus \mathbb{C})$ obviously lands in $\mathrm{SL}(V \oplus \mathbb{C})$.

The following constructions will depend on $\rho'$ (see, e.g., Example 2.1.5), but we won't refer to it explicitly. Define $\rho \colon G \longrightarrow \mathrm{SL}(V) \hookrightarrow \mathrm{GL}(V)$.

If $\mathcal{E}$ is a principal $G$-bundle, the principal $\mathrm{GL}(V)$-bundle associated to $\mathcal{E}$ via $\rho$ is denoted by $\rho_* \mathcal{E}$. Recall that any principal $G$-bundle can be constructed from a pair $(\mathcal{F}, \sigma)$, consisting of a principal $\mathrm{GL}(V)$-bundle $\mathcal{F}$ and a section $\sigma\colon X \longrightarrow \mathcal{F}/G$, by means of the pullback diagram

$$
\begin{array}{ccc}
\sigma^*\mathcal{F} & \longrightarrow & \mathcal{F} \\
\downarrow & & \downarrow \\
X & \xrightarrow{\ \sigma\ } & \mathcal{F}/G.
\end{array}
$$

Moreover, there is a natural isomorphism $\rho_*\sigma^*\mathcal{F} \longrightarrow \mathcal{F}$, and, if we look at $\sigma^*\mathcal{F} \hookrightarrow \rho_*\sigma^*\mathcal{F} \cong \mathcal{F}$ and take the $G$-quotient, we get a section $X = \sigma^*\mathcal{F}/G \longrightarrow \mathcal{F}/G$ which is just $\sigma$ ([13], 4.10).

Suppose we are given $(\mathcal{F}_1, \sigma_1)$ and $(\mathcal{F}_2, \sigma_2)$ and an isomorphism $\varphi\colon \sigma_1^*\mathcal{F}_1 \longrightarrow \sigma_2^*\mathcal{F}_2$, then there is a commutative diagram

$$
\begin{array}{ccccc}
\sigma_1^*\mathcal{F}_1 & \longrightarrow & \rho_*\sigma_1^*\mathcal{F}_1 & \xrightarrow{\cong} & \mathcal{F}_1 \\
\varphi\downarrow & & \rho_*\varphi\downarrow & & \downarrow\psi \\
\sigma_2^*\mathcal{F}_2 & \longrightarrow & \rho_*\sigma_2^*\mathcal{F}_2 & \xrightarrow{\cong} & \mathcal{F}_2.
\end{array}
$$

Taking $G$-quotients, we find $\sigma_2 = \overline{\psi} \circ \sigma_1$, where $\overline{\psi}\colon \mathcal{F}_1/G \longrightarrow \mathcal{F}_2/G$ is the induced isomorphism. Thus, we can identify the set of isomorphism classes of principal $G$-bundles with the set of equivalence classes of pairs $(\mathcal{F}, \sigma)$ where $(\mathcal{F}_1, \sigma_1)$ and $(\mathcal{F}_2, \sigma_2)$ are said to be *equivalent*, if there is an isomorphism $\psi\colon \mathcal{F}_1 \longrightarrow \mathcal{F}_2$ with $\sigma_2 = \overline{\psi} \circ \sigma_1$, $\overline{\psi}\colon \mathcal{F}_1/G \longrightarrow \mathcal{F}_2/G$ being the induced isomorphism.

Let $\mathcal{E}$ be a principal $G$-bundle. The associated $\mathrm{GL}(V)$-bundle $\mathcal{F}$ can be described as follows: Let $\mathcal{V}$ be the vector bundle with fibre $V$ assigned to $\mathcal{E}$ via $\rho$. Then,

$$
\rho_*\mathcal{E} \quad = \quad \underline{\mathrm{Isom}}(V \otimes \mathcal{O}_X, \mathcal{V}) \quad = \quad \bigcup_{x\in X} \mathrm{Isom}(V, \mathcal{V}_x).
$$

Thus, we have a natural inclusion

$$
\rho_*\mathcal{E} \quad \subset \quad \underline{\mathrm{Hom}}(V \otimes \mathcal{O}_X, \mathcal{V}) \quad = \quad \underline{\mathrm{Spec}}\big(\mathrm{Sym}^*(V \otimes \mathcal{V}^\vee)\big).
$$

Now, $V \otimes \mathcal{V}^\vee$ has a natural $\mathrm{GL}(V)$-linearization, inducing the natural $\mathrm{GL}(V)$-action on $\underline{\mathrm{Hom}}(V \otimes \mathcal{O}_X, \mathcal{V})$, and thus, via $\rho$, a $G$-action.

We can form the GIT-quotient

$$
\underline{\mathrm{Hom}}(V \otimes \mathcal{O}_X, \mathcal{V})/\!\!/G \quad = \quad \underline{\mathrm{Spec}}\big(\mathrm{Sym}^*(V \otimes \mathcal{V}^\vee)^G\big).
$$

We must show that $\underline{\mathrm{Hom}}(V \otimes \mathcal{O}_X, \mathcal{V})/\!\!/G$ contains $\rho_*\mathcal{E}/G$ as an open subscheme. Since the formation of the GIT quotient is compatible with closed embeddings, it suffices to look at $\mathrm{Hom}(V, \mathcal{V}_x)$ for $x \in X$. We have to show that the points in $\mathrm{Isom}(V, \mathcal{V}_x)$ are stable. Fix bases for $V$ and $\mathcal{V}_x$, so that we get isomorphisms $\mathbb{C} \cong \bigwedge^{\dim V} V$ and $\mathbb{C} \cong \bigwedge^{\dim V} \mathcal{V}_x$. Then, $\det\colon \mathrm{Hom}(V, \mathcal{V}_x) \longrightarrow \mathbb{C}$, $f \longmapsto \det(f)$, becomes a $G$-invariant function, because $G \subset \mathrm{SL}(V)$. Since

$$
\mathrm{Isom}(V, \mathcal{V}_x) \quad = \quad \big\{\, f \in \mathrm{Hom}(V, \mathcal{V}_x) \mid \det(f) \neq 0 \,\big\},
$$

it suffices to show that the action of $G$ on $\mathrm{Isom}(V, \mathcal{V}_x)$ is closed and all stabilizers are finite, but this is evident.

Finally, the datum of a section $X \longrightarrow \underline{\mathrm{Hom}}(V \otimes \mathcal{O}_X, \mathcal{V}_x)/\!\!/G$ is equivalent to the datum of an $\mathcal{O}_X$-algebra homomorphism $\tau\colon \mathrm{Sym}^*(V \otimes \mathcal{V}^\vee)^G \longrightarrow \mathcal{O}_X$.

*Remark* 2.1.2. In many important special cases, one has $\mathcal{V} \cong \mathcal{V}^\vee$. E.g., for $\mathrm{SO}(r)$ and $\mathrm{Sp}(2r)$ with their standard realizations in $\mathrm{SL}(r)$ and $\mathrm{SL}(2r)$, respectively, because then the associated vector bundle $\mathcal{V}$ of a principal bundle $\mathcal{F}$ comes with an everywhere non-degenerate pairing $\mathcal{V} \otimes \mathcal{V} \longrightarrow \mathcal{O}_X$, providing an isomorphism $\mathcal{V} \longrightarrow \mathcal{V}^\vee$. The same applies to the adjoint representation of any semisimple group. This time, the non-degenerate pairing stems from the Killing form on the Lie algebra bundle.

Set $\mathcal{A} := \mathcal{V}^\vee$. The datum of a pair $(\mathcal{A}, \tau)$, composed of a vector bundle $\mathcal{A}$ with trivial determinant and fibre $V$ and a homomorphism $\tau \colon \mathrm{Sym}^*(V \otimes \mathcal{A})^G \longrightarrow \mathcal{O}_X$ of $\mathcal{O}_X$-algebras can be generalized to "singular" objects. For this, let $\mathcal{A}$ be a torsion free sheaf of rank dim $V$ with trivial determinant. Then,

$$\underline{\mathrm{Hom}}(\mathcal{A}, V^\vee \otimes \mathcal{O}_X) \quad := \quad \underline{\mathrm{Spec}}\big(\mathrm{Sym}^*(V \otimes \mathcal{A})\big)$$

is again a linear space over $X$, carrying a natural $G$-action. The fibre over $x \in X$ of that space is just $\mathrm{Hom}(\mathcal{A}_x, V^\vee) = \mathrm{Hom}(V, \mathcal{A}_x^\vee)$ with the $G$-action being induced via $\rho$ from the natural $\mathrm{GL}(V)$-action. As before, we can form

$$\underline{\mathrm{Hom}}(\mathcal{A}, V^\vee \otimes \mathcal{O}_X)/\!\!/G \quad = \quad \underline{\mathrm{Spec}}\big(\mathrm{Sym}^*(V \otimes \mathcal{A})\big)^G.$$

Any isomorphism $\psi \colon \mathcal{A}_1 \longrightarrow \mathcal{A}_2$ induces an isomorphism $\mathrm{Sym}^* \psi \colon \mathrm{Sym}^*(V \otimes \mathcal{A}_1) \longrightarrow \mathrm{Sym}^*(V \otimes \mathcal{A}_2)$ of $G$-linearized $\mathcal{O}_X$-algebras, yielding an isomorphism

$$\overline{\psi} \colon \mathrm{Sym}^*(V \otimes \mathcal{A}_1)^G \quad \longrightarrow \quad \mathrm{Sym}^*(V \otimes \mathcal{A}_2)^G.$$

Therefore, we look at pairs $(\mathcal{A}, \tau)$ where $\mathcal{A}$ is a torsion free sheaf with trivial determinant of rank dim $V$ and $\tau \colon \mathrm{Sym}^*(V \otimes \mathcal{A})^G \longrightarrow \mathcal{O}_X$ is a non-constant homomorphism of $\mathcal{O}_X$-algebras (that is, $\tau$ is not just the projection onto $\mathcal{O}_X$, or, equivalently, the corresponding section $\sigma \colon X \longrightarrow \underline{\mathrm{Hom}}(\mathcal{A}, V^\vee \otimes \mathcal{O}_X)/\!\!/G$ is not the zero section), and we say that $(\mathcal{A}_1, \tau_1)$ is *equivalent* to $(\mathcal{A}_2, \tau_2)$ if there is an isomorphism $\psi \colon \mathcal{A}_1 \longrightarrow \mathcal{A}_2$ with $\tau_1 = \tau_2 \circ \overline{\psi}$, $\overline{\psi}$ as above. We call $(\mathcal{A}, \tau)$ a *singular principal $G$-bundle*. The *Hilbert polynomial of* $(\mathcal{A}, \tau)$ is just $P(\mathcal{A})$.

*Remark* 2.1.3. Let $(\mathcal{A}, \tau)$ be as above. Then, $\tau$ provides us with a section $\sigma \colon X \longrightarrow \underline{\mathrm{Hom}}(\mathcal{A}, V^\vee \otimes \mathcal{O}_X)/\!\!/G$. Take the fibre product

$$
\begin{array}{ccc}
\mathcal{P} & \longrightarrow & \underline{\mathrm{Hom}}(\mathcal{A}, V^\vee \otimes \mathcal{O}_X) \\
\downarrow & & \downarrow \\
X & \xrightarrow{\ \sigma\ } & \underline{\mathrm{Hom}}(\mathcal{A}, V^\vee \otimes \mathcal{O}_X)/\!\!/G.
\end{array}
$$

Thus, $\mathcal{P}$ is the geometric version of a singular principal $G$-bundle. Let $U$ be the open set over which $\mathcal{A}$ is locally free. Look at $\pi_U \colon H := \underline{\mathrm{Hom}}(\mathcal{A}_{|U}, V^\vee \otimes \mathcal{O}_U) \longrightarrow U$. On $H$, there is the universal homomorphism $\mathfrak{u} \colon \pi_U^* \mathcal{A}_{|U} \longrightarrow V^\vee \otimes \mathcal{O}_H$. The homomorphism $\det(\mathfrak{u}) \colon \mathcal{O}_H \cong \det(\pi_U^* \mathcal{A}_{|U}) \longrightarrow \bigwedge^{\dim V} V^\vee \otimes \mathcal{O}_H \cong \mathcal{O}_H$ is $G$-invariant and, therefore, descends to the quotient $\underline{\mathrm{Hom}}(\mathcal{A}_{|U}, V^\vee \otimes \mathcal{O}_U)/\!\!/G$. Let $\sigma_U \colon X \longrightarrow \underline{\mathrm{Hom}}(\mathcal{A}_{|U}, V^\vee \otimes \mathcal{O}_U)/\!\!/G$ be the section induced by $\tau$. The pullback of the descendant of $\det(\mathfrak{u})$ to $U$ yields a homomorphism $\det \mathcal{A}_{|U} \longrightarrow \mathcal{O}_U$ which can be extended to $\mathfrak{d} \colon \det \mathcal{A} \longrightarrow \mathcal{O}_X$, because $X$ is normal and $\mathrm{codim}(X \setminus U, X) \geq 2$. Therefore, if $\mathfrak{d}$ isn't identically zero, it is an isomorphism, because $\det \mathcal{A}$ is trivial.

**Corollary 2.1.4.** *Either $\mathfrak{d}$ is zero or $\mathcal{P}$ carries the structure of a principal $G$-bundle over the open set $U$ where $\mathcal{A}$ is locally free.*

The latter case looks like a very reasonable generalization of the concept of a torsion free sheaf to the setting of $G$-bundles. But it is not a priori clear that objects of the former type do not have to be included in the compactification of the moduli space of stable principal bundles. In the examples, we will check that, in some special cases, this is not the case. Let us call singular principal $G$-bundles of the latter type *honest singular principal $G$-bundles* in the sequel.

Conversely, given a torsion free sheaf $\mathcal{A}$ with trivial determinant of rank dim $V$ and a section $\sigma_U \colon U \longrightarrow \underline{\mathrm{Isom}}(\mathcal{A}_{|U}, V^\vee \otimes \mathcal{O}_X)/G$ over the open set $U$ where $\mathcal{A}$ is locally free, then $\sigma_U$ extends uniquely to a section $\sigma \colon X \longrightarrow \underline{\mathrm{Hom}}(\mathcal{A}, V^\vee \otimes \mathcal{O}_X) /\!\!/ G$. This is because $\underline{\mathrm{Hom}}(\mathcal{A}, V^\vee \otimes \mathcal{O}_X) /\!\!/ G$ can be embedded as a closed subscheme into a vector bundle.

As an anonymous referee told me, these observations do not suffice to show that the set of isomorphism classes of honest singular principal bundles $(\mathcal{A}, \tau)$ does not depend on the choice of the representation. This is because there might be different extensions $\mathcal{A}$ for $\mathcal{A}_{|U}$ (see the following example which is also due to that referee).

*Example* 2.1.5. First start with the standard representation of $\mathrm{SL}(2)$ on $\mathbb{C}^2$ and look at singular principal bundles $(\mathcal{A}, \tau)$ with $c_2(\mathcal{A}) = 1$ on $\mathbb{P}_2$, such that $\mathcal{A}_{|U} \cong \mathcal{O}_U^{\oplus 2}$, $U := \mathbb{P}_2 \setminus \{\mathrm{pt}\}$. This forces $\mathcal{A} \cong \mathcal{O}_{\mathbb{P}_2} \oplus \mathcal{I}_{\mathrm{pt}}$. If we work with the representation of $\mathrm{SL}(2)$ on $\mathbb{C}^{2 \otimes 2}$ instead, we have to find extensions of $\mathcal{O}_U^{\oplus 4}$ to a torsion free sheaf $\mathcal{B}$ with $c_2(\mathcal{B}) = 2$, and here, $\mathcal{O}_{\mathbb{P}_2}^{\oplus 2} \oplus \mathcal{I}_{\mathrm{pt}}^{\oplus 2}$ and $\mathcal{O}_{\mathbb{P}_2}^{\oplus 3} \oplus \mathcal{I}_Z$ with $Z$ a subscheme of length 2 supported at pt are two different, non-isomorphic choices.

2.2. **Semistability.** We look at the representation $R \colon G \times \mathrm{GL}(r) \longrightarrow \mathrm{GL}(V \otimes \mathbb{C}^r)$, $r := \dim V$, given by $(R(g, g'))(v \otimes w) := \rho(g)(v) \otimes g' \cdot w$. This yields a rational representation of $\mathrm{GL}(r)$ on the algebra $\mathrm{Sym}^*(V \otimes \mathbb{C}^r)^G$, respecting the grading. To see this, let $S^d$ be the homogenous part of degree $d$ in $\mathrm{Sym}^*(V \otimes \mathbb{C}^r)$. Then, we have a representation of $G \times \mathrm{GL}(r)$ on $S^d$. For $(g_1, e)$, $(e, g_2)$ with $g_1 \in G$ and $g_2 \in \mathrm{GL}(r)$ and $v \in S^d$, we obviously have $((g_1, e) \cdot (e, g_2)) \cdot v = (g_1, g_2) \cdot v = ((e, g_2) \cdot (g_1, e)) \cdot v$. This implies that the subspace $(S^d)^G$ of $G$-invariant elements is in fact a $\mathrm{GL}(r)$-module, so that we find the desired rational representation of $\mathrm{GL}(r)$ on the algebra $\mathrm{Sym}^*(V \otimes \mathbb{C}^r)^G$.

Suppose the vector space $\bigoplus_{i=1}^s \mathrm{Sym}^i(V \otimes \mathbb{C}^r)^G$ contains a set of generators for the algebra $\mathrm{Sym}^*(V \otimes \mathbb{C}^r)^G$. Then, we have a representation

$$\mathrm{GL}(r) \longrightarrow \mathrm{GL}(\bigoplus_{i=1}^s \mathrm{Sym}^i(V \otimes \mathbb{C}^r)^G).$$

This representation is not homogeneous in the sense that the restriction to the center of $\mathrm{GL}(r)$ is not of the form $z \cdot \mathbb{E}_r \mapsto z^\alpha \cdot \mathrm{id}$. Therefore, we pass to the induced homogeneous representation

$$t(s) \colon \mathrm{GL}(r) \longrightarrow \mathrm{GL}(\mathbb{U}(s)) \quad \text{with} \quad \mathbb{U}(s) := \bigoplus_{\substack{\underline{d} = (d_1, \ldots, d_s): \\ d_i \geq 0, \ \sum i d_j = s!}} \mathbb{S}^{\underline{d}}$$

and

$$\mathbb{S}^{\underline{d}} \quad := \quad \bigotimes_{i=1}^s \mathrm{Sym}^{d_i}\left(\left(\mathrm{Sym}^i(V \otimes \mathbb{C}^r)\right)^G\right).$$

The upshot is that $\mathbb{P}(V \otimes \mathbb{C}^r)/\!\!/ G$ gets $GL(r)$-equivariantly embedded into the projective space $\mathbb{P}(\mathbb{U}(s))$. Indeed, $\mathbb{P}(V \otimes \mathbb{C}^r)/\!\!/ G = \mathrm{Proj}\big(\mathrm{Sym}^*(V \otimes \mathbb{C}^r)^G\big)$. Let $S^{(s!)}$ be the subalgebra of elements the degree of which is a multiple of $s!$. Then, $\mathrm{Proj}(S^{(s!)}) = \mathrm{Proj}(\mathrm{Sym}^*(V \otimes \mathbb{C}^r)^G)$. We have constructed a degree preserving homomorphism $\mathrm{Sym}^*(\mathbb{U}(s)) \longrightarrow S^{(s!)}$. This map is surjective when $s$ becomes large enough. To see this, fix generators $x_1, ..., x_n$ for $\mathrm{Sym}^*(V \otimes \mathbb{C}^r)^G$ of degrees $d_1, ..., d_n$ and set $d' := \mathrm{lcm}(d_1, ..., d_n)$. Suppose $s$ is so large that $s! \geq n \cdot d'$. Now, add generators $x_{n+1}, ..., x_k$, $k := s!/d'$, of degree $d_n$. From this point, continue as [12], III.8, Proof of Lemma, p. 282. We obtain a closed embedding

$$\mathbb{P}(V \otimes \mathbb{C}^r)/\!\!/ G = \mathrm{Proj}(S^{(s!)}) \hookrightarrow \mathrm{Proj}\big(\mathrm{Sym}^*(\mathbb{U}(s))\big) = \mathbb{P}(\mathbb{U}(s)).$$

Now, we are in a situation very similar to the one considered in [15] and [4].

Let $(\mathcal{A}, \tau)$ be a singular principal bundle. Note that we require $\tau$ to be non-constant. Let $\sigma \colon X \longrightarrow \mathrm{Hom}(\mathcal{A}, V^\vee \otimes \mathcal{O}_X)/\!\!/ G$ be the associated section. Let $U$ be a non-empty open set over which $\mathcal{A}$ is trivial and the section $\sigma$ is non-zero. Then, over $U$, we have the map

$$\overline{\sigma}_U(s) \colon U \longrightarrow \mathbb{P}\big(V \otimes \mathcal{A}_{|U}\big)/\!\!/ G \hookrightarrow \mathbb{P}(\mathbb{U}(s)) \times U \longrightarrow \mathbb{P}(\mathbb{U}(s)).$$

A *weighted filtration of* $\mathcal{A}$ is a pair $(\mathcal{A}^\bullet, \underline{\alpha})$ composed of a filtration $0 \subset \mathcal{A}^1 \subset \cdots \subset \mathcal{A}^u \subset \mathcal{A}$ by saturated subsheaves, i.e., the quotients $\mathcal{A}/\mathcal{A}^i$ are again torsion free, and a vector $\underline{\alpha} = (\alpha_1, ..., \alpha_u)$ of positive rational numbers. Set

$$M(\mathcal{A}^\bullet, \underline{\alpha}) \quad := \quad \sum_{j=1}^u \alpha_j \big(P(\mathcal{A}) \operatorname{rk} \mathcal{A}^j - P(\mathcal{A}^j) \operatorname{rk} \mathcal{A}\big).$$

Making $U$ slightly smaller, we get a filtration $0 \subset \mathcal{A}^1_{|U} \subset \cdots \subset \mathcal{A}^u_{|U} \subset \mathcal{A}_{|U}$ of $\mathcal{A}$ by subbundles. Then, we can define $\mu\big(\mathcal{A}^\bullet_{|U}, \underline{\alpha}; \overline{\sigma}_U(s)\big)$ as in [15], page 18, and

$$\mu\big(\mathcal{A}^\bullet, \underline{\alpha}; \tau\big) \quad := \quad \frac{1}{s!} \mu\big(\mathcal{A}^\bullet_{|U}, \underline{\alpha}; \overline{\sigma}_U(s)\big).$$

*Remark* 2.2.1. As in [15], the value of $\mu\big(\mathcal{A}^\bullet_{|U}, \underline{\alpha}; \overline{\sigma}_U(s)\big)$ doesn't depend on the choices involved. The factor $(1/s!)$ has been thrown in to make the definition independent of the chosen $s$. In fact, $\overline{\sigma}_U(s+1)$ is $\overline{\sigma}_U(s)$ followed by the $(s+1)$-st Veronese embedding $\mathbb{P}(\mathbb{U}(s)) \hookrightarrow \mathbb{P}(\mathbb{U}(s+1))$.

Fix a positive polynomial $\delta \in \mathbb{Q}[x]$ of degree at most $\dim X - 1$. A singular principal $G$-bundle $(\mathcal{A}, \tau)$ will be called *δ-(semi)stable*, if for every weighted filtration $(\mathcal{A}^\bullet, \underline{\alpha})$ of $\mathcal{A}$, the inequality

$$M(\mathcal{A}^\bullet, \underline{\alpha}) + \delta \cdot \mu\big(\mathcal{A}^\bullet, \underline{\alpha}; \tau\big) \quad (\succeq) \quad 0$$

is satisfied.

*Remark* 2.2.2. i) Fix a polynomial $P$, and let $(\mathcal{A}, \tau)$ be a $\delta$-semistable singular $G$-principal bundle with Hilbert polynomial $P$. For any saturated subsheaf $0 \subsetneq \mathcal{A}' \subsetneq \mathcal{A}$ of $\mathcal{A}$, we look at the weighted filtration $(0 \subsetneq \mathcal{A}' \subsetneq \mathcal{A}, (1))$. As $\mu\big(0 \subsetneq \mathcal{A}' \subsetneq \mathcal{A}, (1); \overline{\sigma}_U(s)\big) \leq s!(r-1)$ ([15], Lemma 1.2.6), we see

$$\mu\big(0 \subsetneq \mathcal{A}' \subsetneq \mathcal{A}, (1); \tau\big) \quad \leq \quad r-1.$$

This implies

$$\mu_{\max}(\mathcal{A}) \quad \leq \quad \mu(\mathcal{A}) + \overline{\delta} \cdot (r-1),$$

where $\overline{\delta}$ is the coefficient of $x^{\dim X-1}$ in $\delta$. Therefore, the set of isomorphism classes of torsion free sheaves $\mathcal{A}$, belonging to $\delta$-semistable singular principal bundles with Hilbert polynomial $P$, is bounded, by Maruyama's theorem [9].

ii) Ramanathan and Subramanian have defined semi- and quasi-stability for principal $G$-bundles over $X$ [14]. Moreover, they have shown that $\rho_*\mathcal{E}$ is a Mumford-semistable (Mumford-quasi-stable) vector bundle, if $\mathcal{E}$ is a semistable (quasi-stable) principal $G$-bundle ([14], Thm. 3). Here, a vector bundle is called *Mumford-quasi-stable*, if it is a direct sum of Mumford-stable bundles.

Based on this, in the case of the adjoint representation, Hyeon has defined a principal $G$-bundle to be stable, if the associated vector bundle is stable.

To see that our definition of $\delta$-semistability extends Hyeon's definition, let $\mathcal{A}$ be a Gieseker-semistable torsion free sheaf, e.g., a Mumford-stable vector bundle. Let $(\mathcal{A}, \tau)$ be an honest singular $G$-bundle. We assert that $\mu(\mathcal{A}^\bullet, \underline{\alpha}; \tau) \geq 0$ for every weighted filtration $(\mathcal{A}^\bullet, \underline{\alpha})$ of $\mathcal{A}$. For this it suffices to show that the images of isomorphisms in $\mathbb{P}(V \otimes \mathbb{C}^r)/\!/G$ are $\mathrm{SL}(r)$-semistable. This amounts to the fact that isomorphisms in $\mathbb{P}(V \otimes \mathbb{C}^r)$ are $(G \times \mathrm{SL}(r))$-semistable. We claim that they are even $(\mathrm{SL}(V) \times \mathrm{SL}(r))$-semistable. Fix a basis for $V$. Then, we have to check that isomorphisms in $\mathbb{P}\big(\mathrm{Hom}(\mathbb{C}^r, \mathbb{C}^{r^\vee})^\vee\big)$ are semistable w.r.t. the $(\mathrm{SL}(r) \times \mathrm{SL}(r))$-action $(g, g') \cdot [f] = [g'^{-1} \circ f \circ g'^{-1}]$. This is true, because the determinant does not vanish on them. We conclude that $(\mathcal{A}, \tau)$ will be $\delta$-semistable for every polynomial $\delta$. Moreover,

$$M(\mathcal{A}^\bullet, \underline{\alpha}) + \delta \cdot \mu(\mathcal{A}^\bullet, \underline{\alpha}; \tau) \quad = \quad 0$$

can only occur if both $M(\mathcal{A}^\bullet, \underline{\alpha})$ and $\mu(\mathcal{A}^\bullet, \underline{\alpha}; \tau)$ are zero.

There is also a kind of a converse to this observation:

**Lemma 2.2.3.** i) *There is a constant polynomial $\delta_{\mathrm{Gies}}$, such that, for every polynomial $\delta' \preceq \delta_{\mathrm{Gies}}$ and every $\delta'$-semistable singular principal bundle $(\mathcal{A}, \tau)$, the sheaf $\mathcal{A}$ is itself Gieseker-semistable.*

ii) *There is a polynomial $\delta_\mu$ of degree exactly $\dim X - 1$, such that, for every polynomial $\delta' \preceq \delta_\mu$ and every $\delta'$-semistable singular principal bundle $(\mathcal{A}, \tau)$, the sheaf $\mathcal{A}$ is itself Mumford-semistable.*

*Proof.* We show ii), i) being similar. Let $(\mathcal{A}, \tau)$ be a singular principal bundle. Suppose $\mathcal{A}' \subsetneq \mathcal{A}$ is a non-trivial saturated subsheaf with $\mu(\mathcal{A}') > \mu(\mathcal{A})$. Take the weighted filtration $(0 \subsetneq \mathcal{A}' \subsetneq \mathcal{A}, (1))$. Then, $M(0 \subsetneq \mathcal{A}' \subsetneq \mathcal{A}, (1))$ is a negative polynomial over the integers of degree exactly $\dim X - 1$. Since $\mu(0 \subsetneq \mathcal{A}' \subsetneq \mathcal{A}, (1); \tau) \leq r - 1$, the assertion is obvious. $\qquad\square$

iii) Consider, for $\underline{d} = (d_1, ..., d_s)$ with $\sum i d_i = s!$, the natural homomorphism

$$\bigotimes_{i=1}^s (V \otimes \mathcal{A})^{\otimes d_i i} \longrightarrow \bigotimes_{i=1}^s \mathrm{Sym}^{d_i}\big(\mathrm{Sym}^i(V \otimes \mathcal{A})\big) \longrightarrow \mathbb{S}^{\underline{d}} \longrightarrow \mathcal{O}_X.$$

These add to

$$\varphi_\tau\colon \Big((V \otimes \mathcal{A})^{\otimes s!}\Big)^{\oplus N} \longrightarrow \mathbb{U}(s) \longrightarrow \mathcal{O}_X.$$

Here, the second homomorphism comes from the Reynolds operator and the third one is induced by $\tau$. Now, $(\mathcal{A}, \varphi_\tau)$ is just a tensor field in the sense of [4]. Then, $(\mathcal{A}, \tau)$is $\delta$-(semi)stable, if and only if $(\mathcal{A}, \varphi_\tau)$ is a $(\delta/s!)$-(semi)stable tensor field as defined in [4]. This is because the quantity $-\mu(\varphi_\tau, \mathcal{A}^\bullet, \underline{\alpha})$ defined in [4] equals $\mu\big(\mathcal{A}_{|U}^{\bullet}, \underline{\alpha}; \widetilde{\sigma}_U(s)\big)$ as defined in [15] w.r.t.

$$\widetilde{\sigma}_U(s)\colon U \xrightarrow{\overline{\sigma}_U(s)} \mathbb{P}(\mathbb{U}(s)) \hookrightarrow \mathbb{P}\Big(\big((V \otimes \mathcal{A})^{\otimes s!}\big)^{\oplus N}\Big)$$

and, obviously,

$$\mu\big(\mathcal{A}_{|U}^{\bullet}, \underline{\alpha}; \widetilde{\sigma}_U(s)\big) \quad = \quad \mu\big(\mathcal{A}_{|U}^{\bullet}, \underline{\alpha}; \overline{\sigma}_U(s)\big).$$

2.3. **Families and the moduli functors.** Let $S$ be a noetherian scheme. A *family of singular G-principal bundles parametrized by* $S$ is a pair $(\mathcal{A}_S, \tau_S)$ consisting of an $S$-flat family $\mathcal{A}_S$ of torsion free coherent sheaves on $X$ and a homomorphism $\tau_S\colon \mathrm{Sym}^*(V \otimes \mathcal{A}_S)^G \longrightarrow \mathcal{O}_{S \times X}$. The definition of *equivalence of families* is left to the reader.

*Remark* 2.3.1. i) By the base change properties discussed in Section 1.4, there is a pullback operation on families of singular principal bundles via base change morphisms $T \longrightarrow S$.

ii) The algebra $\mathrm{Sym}^*(V \otimes \mathcal{A}_S)$ is naturally graded with $\mathrm{Sym}^d(V \otimes \mathcal{A}_S)$ as the homogeneous part of degree $d$. The algebra $\mathrm{Sym}^*(V \otimes \mathcal{A}_S)^G$ inherits this grading, i.e.,

$$\mathrm{Sym}^*(V \otimes \mathcal{A}_S)^G \quad = \quad \bigoplus_{d \geq 0} \mathrm{Sym}^d(V \otimes \mathcal{A}_S)^G.$$

Therefore, $\tau$ breaks up into a collection $\tau_d\colon \mathrm{Sym}^d(V \otimes \mathcal{A}_S)^G \longrightarrow \mathcal{O}_{S \times X}$, $d \geq 0$, of maps between coherent $\mathcal{O}_{S \times X}$-modules. Moreover, since $\big(\mathrm{Sym}^*(V \otimes \mathcal{A}_S)\big)^G$ is finitely generated over $\mathcal{O}_{S \times X}$, finitely many of the $\tau_d$ will suffice to reconstruct $\tau$. As I was informed by a referee, the module $\mathrm{Sym}^d(V \otimes \mathcal{A}_S)^G$ has in general no chance of being $S$-flat. Fortunately, Proposition 1.5.1 of Gómez and Sols shows that this is not a problem. This is the basic observation to start the "standard procedure of constructing moduli spaces with GIT".

With these observations, we can introduce the moduli functors for $\delta$-(semi)stable singular principal $G$-bundles with Hilbert polynomial $P$:

$$\underline{\mathbf{M}}(\rho)_P^{\delta-(s)s}\colon \quad \underline{\mathrm{Sch}}_{\mathbb{C}} \quad \longrightarrow \quad \underline{\mathrm{Set}}$$
$$S \quad \longmapsto \quad \left\{ \begin{array}{c} \text{Equivalence classes of families} \\ \text{of } \delta\text{-(semi)stable singular principal} \\ G\text{-bundles parametrized by } S \end{array} \right\}.$$

## 3. THE MAIN RESULT

**Theorem.** *There exist a projective scheme* $\mathcal{M}(\rho)_P^{\delta-ss}$ *and an open subscheme* $\mathcal{M}(\rho)_P^{\delta-s} \subset \mathcal{M}(\rho)_P^{\delta-ss}$ *together with natural transformations*

$$\vartheta^{(s)s}\colon \underline{\mathbf{M}}(\rho)_P^{\delta-(s)s} \longrightarrow h_{\mathcal{M}(\rho)_P^{\delta-(s)s}},$$

*satisfying the following universal properties:*

1. *For every scheme* $\mathcal{N}$ *and every natural transformation* $\vartheta'\colon \underline{\mathbf{M}}(\rho)_P^{\delta-(s)s} \longrightarrow h_{\mathcal{N}}$, *there exists a unique morphism* $\varphi\colon \mathcal{M}(\rho)_P^{\delta-(s)s} \longrightarrow \mathcal{N}$ *with* $\vartheta' = h(\varphi) \circ \vartheta^{(s)s}$.
2. *The scheme* $\mathcal{M}(\rho)_P^{\delta-s}$ *is a coarse moduli space for the functor* $\underline{\mathbf{M}}(\rho)_P^{\delta-s}$.

## 4. Proof of the main result

Let $\mathfrak{A}$ be the bounded set of isomorphism classes of torsion free sheaves with fixed Hilbert polynomial $P$ and trivial determinant occuring in $\delta$-semistable singular principal bundles (see Rem. 2.2.2 i)).

4.1. **The parameter space.** By the boundedness of $\mathfrak{A}$, there exists an $n_0$, such that for every $n \geq n_0$ and every torsion free sheaf $\mathcal{A}$ with $[\mathcal{A}] \in \mathfrak{A}$

- $H^i(\mathcal{A}(n)) = 0, i > 0$,
- $\mathcal{A}(n)$ is generated by global sections.

Fix such an $n$, and let $W$ be a complex vector space of dimension $P(n)$. Let $\mathfrak{Q}$ be the quot scheme, parametrizing quotients $W \otimes \mathcal{O}_X(-n) \longrightarrow \mathcal{A}$ where $\mathcal{A}$ is a coherent sheaf on $X$ with Hilbert polynomial $P$ and trivial determinant. Denote the universal quotient by

$$\mathfrak{q}_{\mathfrak{Q}} \colon W \otimes \pi_X^* \mathcal{O}_X(-n) \longrightarrow \mathcal{A}_{\mathfrak{Q}}.$$

Look at the following commutative diagram

$$\begin{array}{ccc}
\operatorname{Sym}^*\big(V \otimes W \otimes \pi_X^* \mathcal{O}_X(-n)\big) & \longrightarrow & \operatorname{Sym}^*\big(V \otimes \mathcal{A}_{\mathfrak{Q}}\big) \\
\downarrow & & \downarrow \\
\operatorname{Sym}^*\big(V \otimes W \otimes \pi_X^* \mathcal{O}_X(-n)\big)^G & \longrightarrow & \operatorname{Sym}^*\big(V \otimes \mathcal{A}_{\mathfrak{Q}}\big)^G
\end{array}$$

where the vertical maps come from the Reynolds operator. Thus, we have a surjective homomorphism of $\mathcal{O}_{\mathfrak{Q} \times X}$-modules

$$h \colon \operatorname{Sym}^*\big(V \otimes W \otimes \pi_X^* \mathcal{O}_X(-n)\big) \longrightarrow \operatorname{Sym}^*\big(V \otimes \mathcal{A}_{\mathfrak{Q}}\big)^G.$$

For every affine open subset $U = \operatorname{Spec} R \subset \mathfrak{Q} \times X$, the $R$-algebra $\operatorname{Sym}^*\big(V \otimes \mathcal{A}_{\mathfrak{Q}}(U)\big)^G$ is a finitely generated $\mathbb{C}$-algebra. Thus, there is an $s$, such that

$$h\Big(\bigoplus_{i=1}^s \operatorname{Sym}^i\big(V \otimes W \otimes \pi_X^* \mathcal{O}_X(-n)\big)\Big)$$

will contain a set of generators for this algebra. Since $\mathfrak{Q} \times X$ is quasi-compact, we can find an $s$ working for any affine open set.

A homomorphism

$$k \colon \bigoplus_{i=1}^s \operatorname{Sym}^i\big(V \otimes W \otimes \mathcal{O}_X(-n)\big) \longrightarrow \mathcal{O}_X$$

breaks into a family of homomorphisms

$$k^i \colon \operatorname{Sym}^i(V \otimes W) \otimes \mathcal{O}_X \longrightarrow \mathcal{O}_X(in), \quad i = 1, ..., s.$$

These are determined by the maps

$$\kappa^i := H^0(k^i) \colon \operatorname{Sym}^i(V \otimes W) \longrightarrow H^0(\mathcal{O}_X(in)), \quad i = 1, ..., s,$$

on global sections. Therefore, our first approximation of the parameter space is

$$\overline{\mathfrak{Y}} \quad := \quad \mathfrak{Q} \times \bigoplus_{i=1}^s \operatorname{Hom}\Big(\operatorname{Sym}^i(V \otimes W), H^0\big(\mathcal{O}_X(in)\big)\Big).$$

We now have to single out those points $([q], [\kappa])$ where $\kappa$ comes from an algebra homomorphism $\mathrm{Sym}^*(V \otimes \mathcal{A}_{\overline{\mathfrak{Q}}|[q] \times X})^G \longrightarrow \mathcal{O}_X$. In order to do so, we observe that, on $\overline{\mathfrak{Q}} \times X$, there are universal homomorphisms

$$\widetilde{\varphi}^i \colon \mathrm{Sym}^i(V \otimes W) \otimes \mathcal{O}_{\overline{\mathfrak{Q}} \times X} \longrightarrow H^0(\mathcal{O}_X(in)) \otimes \mathcal{O}_{\overline{\mathfrak{Q}} \times X}, \quad i = 1, ..., s.$$

Compose these with the pullbacks of the evaluation maps $H^0(\mathcal{O}_X(in)) \otimes \mathcal{O}_X \longrightarrow \mathcal{O}_X(in)$ to get

$$\varphi^i \colon \mathrm{Sym}^i(V \otimes W) \otimes \mathcal{O}_{\overline{\mathfrak{Q}} \times X} \longrightarrow \pi_X^* \mathcal{O}_X(in), \quad i = 1, ..., s.$$

These can be put together to

$$\varphi \colon \mathcal{V}_{\overline{\mathfrak{Q}}} := \bigoplus_{i=1}^{s} \mathrm{Sym}^i\big(V \otimes W \otimes \pi_X^* \mathcal{O}_X(-n)\big) \longrightarrow \mathcal{O}_{\overline{\mathfrak{Q}} \times X}.$$

Graduate the symmetric algebra $\mathrm{Sym}^*(\mathcal{V}_{\overline{\mathfrak{Q}}})$ by assigning the weight $i$ to the elements in $\mathrm{Sym}^i(...)$. Now, $\varphi$ gives us a homomorphism of $\mathcal{O}_{\overline{\mathfrak{Q}} \times X}$-algebras

$$\widetilde{\tau}_{\overline{\mathfrak{Q}}} \colon \mathrm{Sym}^*(\mathcal{V}_{\overline{\mathfrak{Q}}}) \longrightarrow \mathcal{O}_{\overline{\mathfrak{Q}} \times X}.$$

On the other hand, we have a surjective homomorphism of graded $\mathcal{O}_{\overline{\mathfrak{Q}} \times X}$-algebras

$$\beta \colon \mathrm{Sym}^*(\mathcal{V}_{\overline{\mathfrak{Q}}}) \longrightarrow \mathrm{Sym}^*\big(V \otimes \pi_{\overline{\mathfrak{Q}} \times X}^* \mathcal{A}_{\overline{\mathfrak{Q}}}\big)^G.$$

Now, the subscheme $\mathfrak{Y}$ will be defined by the condition that $\widetilde{\tau}_{\overline{\mathfrak{Q}}}$ vanish on $\ker \beta$. To see that this is a closed subscheme of $\overline{\mathfrak{Q}}$, note that

$$\mathfrak{Y} \quad = \quad \bigcap_{d \geq 0} \mathfrak{Y}^d$$

with

$$\mathfrak{Y}^d \quad := \quad \big\{ y \in \overline{\mathfrak{Q}} \,\big|\, \widetilde{\tau}_{\overline{\mathfrak{Q}}|\{y\} \times X}^d \colon \ker \beta_{|\{y\} \times X}^d \longrightarrow \mathcal{O}_X \text{ is trivial} \big\}.$$

Here, $\widetilde{\tau}_{\overline{\mathfrak{Q}}}^d$ and $\beta^d$ stand for the degree $d$ component of the respective homomorphism. By Proposition 1.5.1, $\mathfrak{Y}^d$ is a closed subscheme of $\overline{\mathfrak{Q}}$, $d \geq 0$, whence so is $\mathfrak{Y}$. Let $\mathcal{A}_{\mathfrak{Y}}$ be the pullback of $\mathcal{A}_{\overline{\mathfrak{Q}}}$ to $\mathfrak{Y} \times X$ and

$$\tau_{\mathfrak{Y}} \colon \mathrm{Sym}^*\big(V \otimes \mathcal{A}_{\mathfrak{Y}}\big)^G \longrightarrow \mathcal{O}_{\mathfrak{Y} \times X}$$

the homomorphism of $\mathcal{O}_{\mathfrak{Y} \times X}$-algebras induced by $\widetilde{\tau}_{\overline{\mathfrak{Q}}|\mathfrak{Y} \times X}$. We call $(\mathcal{A}_{\mathfrak{Y}}, \tau_{\mathfrak{Y}})$ the *universal family of singular principal bundles*. The following is deduced with standard arguments from the universal property of $\mathfrak{Y}$ (see Proposition 1.5.1).

**Proposition 4.1.1** (Local universal property)**.** *Let $S$ be a scheme of finite type over $\mathbb{C}$, and $(\mathcal{A}_S, \tau_S)$ a family of $\delta$-semistable singular principal bundles parametrized by $S$. Then, there exists an open covering $S_i$, $i \in I$, of $S$, and morphisms $\beta_i \colon S_i \longrightarrow \mathfrak{Y}$, $i \in I$, such that the restriction of the family $(\mathcal{A}_S, \tau_S)$ to $S_i \times X$ is equivalent to the pullback of $(\mathcal{A}_{\mathfrak{Y}}, \tau_{\mathfrak{Y}})$ via $\beta_i \times \mathrm{id}_X$, for all $i \in I$.*

4.2. **The group action.** One has a natural $\mathrm{GL}(W)$-action on $\mathfrak{Q}$, and the universal quotient $\mathcal{A}_\mathfrak{Q}$ is linearized w.r.t. this group action. Thus, the algebras $\mathrm{Sym}^*(V \otimes \mathcal{A}_\mathfrak{Q})$ and $\mathrm{Sym}^*(V \otimes \mathcal{A}_\mathfrak{Q})^G$ are also $\mathrm{GL}(W)$-linearized. The natural $\mathrm{GL}(W)$-action on $\overline{\mathfrak{Y}}$, therefore, leaves the parameter space $\mathfrak{Y}$ invariant and induces an action

$$\Gamma\colon \mathrm{GL}(W) \times \mathfrak{Y} \longrightarrow \mathfrak{Y}.$$

Let $\mathfrak{Y}^0$ be the open part of $\mathfrak{Y}$ consisting of the points $([q], [f])$, such that $H^0([q])\colon W \longrightarrow H^0(\mathcal{A}_{\mathfrak{Q}|\{[q]\}\times X}(n))$ is an isomorphism and $\mathcal{A}_{\mathfrak{Q}|\{[q]\}\times X}$ is torsion free. The following is again standard.

**Proposition 4.2.1** (Gluing property). *Let $S$ be a scheme of finite type over $\mathbb{C}$ and $\beta_{1,2}\colon S \longrightarrow \mathfrak{Y}^0$ two morphisms, such that the pullbacks of $(\mathcal{A}_\mathfrak{Y}, \tau_\mathfrak{Y})$ via $\beta_1 \times \mathrm{id}_X$ and $\beta_2 \times \mathrm{id}_X$ are equivalent. Then, there exists a map $\Xi\colon S \longrightarrow \mathrm{GL}(W)$, such that the morphism $\beta_2$ equals the morphism*

$$S \xrightarrow{\Xi \times \beta_1} \mathrm{GL}(W) \times \mathfrak{Y} \xrightarrow{\Gamma} \mathfrak{Y}.$$

We can view the $\mathrm{GL}(W)$-action also as a $(\mathbb{C}^* \times \mathrm{SL}(W))$-action. First, we divide by the $\mathbb{C}^*$-action. The quotient $\mathfrak{Y} /\!/ \mathbb{C}^*$ is a closed subscheme of

$$\mathfrak{Q} \times \left( \left( \bigoplus_{i=1}^s \mathrm{Hom}\big(\mathrm{Sym}^i(V \otimes W), H^0(\mathcal{O}_X(in))\big) \right) /\!/ \mathbb{C}^* \right).$$

The right hand factor embeds into

$$\mathbb{T} \quad := \quad \mathbb{P}\Big( \bigoplus_{\substack{\underline{d}=(d_1,\dots,d_s):\\ d_i \geq 0,\ \sum i d_i = s!}} \mathbb{T}_{\underline{d}}^\vee \Big)$$

with

$$\mathbb{T}_{\underline{d}} \quad := \quad \bigotimes_{i=1}^s \mathrm{Hom}\Big(\mathrm{Sym}^{d_i}\mathrm{Sym}^i(V \otimes W), \mathrm{Sym}^{d_i}H^0(\mathcal{O}_X(in))\Big).$$

The space $\mathbb{T}$ naturally maps to

$$\mathfrak{Z} \quad := \quad \mathbb{P}\Big( \mathrm{Hom}\big( ((V \otimes W)^{\otimes s!})^{\oplus N}, H^0(\mathcal{O}_X(s! \cdot n)) \big)^\vee \Big).$$

The induced morphism

$$\mathfrak{Y} /\!/ \mathbb{C}^* \longrightarrow \mathfrak{Q} \times \mathfrak{Z}$$

is $\mathrm{SL}(W)$-equivariant and injective. In fact, it is the map which associates to $([q\colon W \otimes \mathcal{O}_X(-n) \longrightarrow \mathcal{A}], [\tau])$ the quotient $[q]$ and the class of the map

$$
\begin{array}{cc}
((V \otimes W)^{\otimes s!})^{\oplus N} & \xrightarrow{(\mathrm{id}_V \otimes H^0(q(n))^{\otimes s!})^{\oplus N}} & ((V \otimes H^0(\mathcal{A}(n)))^{\otimes s!})^{\oplus N} \\
& \xrightarrow{H^0(\varphi_\tau(s! \cdot n))} & H^0(\mathcal{O}_X(s! \cdot n))
\end{array}
$$

where $\varphi_\tau$ is the associated tensor field (2.2.2). Note that since we do not have twists by line bundles in our tensor fields and no additional data of coherent sheaves in our setting, the space $\mathfrak{Z}$ is indeed the correct specialization of the space *"P"* in [4], p. 12, and $\mathfrak{Y} /\!/ \mathbb{C}^*$ lands in *"Z'"*.

By Proposition 4.1.1 and 4.2.1, the main result now follows from

**Theorem 4.2.2.** *There is a linearization of the $\mathrm{SL}(W)$-action on $\mathfrak{Q} \times \mathfrak{Z}$, such that for the pullback $\eta$ of this linearization one has:*

1. *All $\eta$-semistable points lie in $\mathfrak{Y}^0$.*

2. *A point $y \in \mathfrak{Y}^0$ is $\eta$-(semi)stable if and only if the restriction of the universal singular principal bundle to $\{y\} \times X$ is $\delta$-(semi)stable.*

*Proof.* We have constructed an injective and proper, whence finite, morphism $\mathfrak{Y} /\!/ \mathbb{C}^* \longrightarrow Z'$, $Z'$ as in [4]. Since we work with the pullback of the linearization from $Z'$ to $\mathfrak{Y} /\!/ \mathbb{C}^*$, the (semi)stable points in $\mathfrak{Y} /\!/ \mathbb{C}^*$ are just the pullbacks of the (semi)stable points in $Z'$ (see [11]). By virtue of the relation with tensor fields established in Remark 2.2.2, this theorem follows from Theorem 3.5 in [4]. □

### 4.3. A word about the closed points.

Recall that the points in the GIT quotient of a projective algebraic scheme $Y$ correspond to the closed orbits in the set of semistable points, and that every point in $\overline{\mathrm{Orb}}(y) \cap Y^{ss}$, with $y \in Y^{ss}$ can be obtained as $\lim_{z\to\infty} \lambda(z) \cdot y$ where $\lambda \colon \mathbb{C}^* \longrightarrow G$ is a one parameter subgroup with $\mu(\lambda, y) = 0$. Suppose $y = ([q \colon W \otimes \mathcal{O}_X(-n) \longrightarrow \mathcal{A}], \tau) \in \mathfrak{Y}^0$. A one parameter subgroup $\lambda \colon \mathbb{C}^* \longrightarrow \mathrm{GL}(W)$ defines a weighted filtration $(\mathcal{A}^\bullet, \underline{\alpha})$ of $\mathcal{A}$ with $M(\mathcal{A}^\bullet, \underline{\alpha}) + \delta \cdot \mu(\mathcal{A}^\bullet, \underline{\alpha}; \tau) = 0$, and every weighted filtration with this property arises in this way. For the point $y_\infty = ([q \colon W \otimes \mathcal{O}_X(-n) \longrightarrow \mathcal{A}_\infty], \tau_\infty) := \lim_{z\to\infty} \lambda(z) \cdot y$ one has $\mathcal{A}_\infty = \mathrm{gr}(\mathcal{A}^\bullet)$, the graded object associated with the filtration $\mathcal{A}^\bullet$, and the associated tensor field $\varphi_{\tau_\infty}$ is described in [4].

Finally, assume that $X$ is a smooth projective curve and $(\mathcal{A}, \tau)$ comes from a quasi-stable principal $G$-bundle $\mathcal{F}$. This implies that $\mathcal{F}$ has a reduction to a Levi component $L$ of a parabolic subgroup of $G$. Since the associated vector bundle is Mumford-quasi-stable, we find $\rho(L) \subset \mathrm{GL}(V_1) \times \cdots \times \mathrm{GL}(V_t) \subset \mathrm{GL}(V)$, and the associated vector bundles $\mathcal{A}_i$ with fibre $V_i$, $i = 1, ..., t$, are all Mumford-stable of degree 0. Translated to our setting, this means that $\tau$ is of the following form

$$\mathrm{Sym}^*(V \otimes \mathcal{A})^G \subset \mathrm{Sym}^*(V \otimes \mathcal{A})^L \longrightarrow \bigoplus_{i=1}^{t} \mathrm{Sym}^*(V_i \otimes \mathcal{A}_i)^L \longrightarrow \mathcal{O}_X.$$

We have explained in Remark 2.2.2 that for any weighted filtration $(\mathcal{A}^\bullet, \underline{\alpha})$ with $M(\mathcal{A}^\bullet, \underline{\alpha}) + \delta \cdot \mu(\mathcal{A}^\bullet, \underline{\alpha}; \tau) = 0$, we must have $M(\mathcal{A}^\bullet, \underline{\alpha}) = 0$. This implies that $\mathcal{A}^\bullet$ is of the form $0 \subset \mathcal{A}_{\iota_1} \subset \mathcal{A}_{\iota_1} \oplus \mathcal{A}_{\iota_2} \subset \cdots \subset \mathcal{A}_{\iota_1} \oplus \cdots \oplus \mathcal{A}_{\iota_u} \subset \mathcal{A}$ where the $\iota_i \in \{1, ..., t\}$ are pairwise distinct. By the above factorization of $\tau$, the associated tensor field $\varphi_\tau$ is non-zero only on the components $\bigotimes_{i=1}^t V_i \otimes \mathcal{A}_i$. This implies that $z := ([q \colon W \otimes \mathcal{O}_X(-n) \longrightarrow \mathcal{A}], [\varphi_\tau]) \in \mathfrak{Q} \times \mathfrak{J}$ is a polystable point, i.e., $\overline{\mathrm{SL}(W) \cdot z} \cap (\mathfrak{Q} \times \mathfrak{J})^{ss} = \mathrm{SL}(W) \cdot z$, whence so is $y$. To summarize, we see that quasi-stable principal $G$-bundles define closed points in $\mathcal{M}(\rho)_P^{\delta-ss}$.

## 5. EXAMPLES

### 5.1. The case of curves.

It is easy to see that there is an open subset $\mathfrak{P}\mathfrak{B} \subset \mathfrak{Y}^0$, consisting of those $([q \colon W \otimes \mathcal{O}_X(-n) \longrightarrow \mathcal{A}], \tau)$ where $\mathcal{A}$ is a locally free, Gieseker-semistable sheaf, and $\tau$ induces a section $\sigma \colon X \longrightarrow \underline{\mathrm{Isom}}(\mathcal{A}, V \otimes \mathcal{O}_X)/G$. Suppose $X$ is a smooth projective curve. Then, we have seen that $\mathfrak{P}\mathfrak{B} \subset \mathfrak{Y}^{ss}$ is a saturated subset, i.e., for every $y \in \mathfrak{P}\mathfrak{B}$, $\overline{\mathrm{GL}(W) \cdot y} \cap \mathfrak{Y}^{ss} \subset \mathfrak{P}\mathfrak{B}$. Therefore, $\mathcal{P}\mathcal{B} := \mathfrak{P}\mathfrak{B} /\!/ \mathrm{GL}(W)$ is an open subset of $\mathcal{M}(\rho)_P^{\delta-ss}$. As it is the moduli space of semistable principal $G$-bundles of "fixed topological type", it equals Ramanathan's moduli space [13]. In particular, it is projective.

5.2. **Classical groups.** If we take $G = SO(r)$ or $G = Sp(2r)$, we can compare our construction with the one of Gómez and Sols [4].

Let us look at the case $G = SO(r) \stackrel{\rho'}{\subset} SL(r)$. Set $V := \mathbb{C}^r$, and let $V \longrightarrow V^\vee$ be the isomorphism coming from the given non-degenerate pairing. Observe that we have an $SO(r)$-invariant morphism $\text{Hom}(\mathbb{C}^r, V^\vee) \longrightarrow \text{Sym}_n \times \mathbb{C}$, $f \longmapsto (f^t \circ f, \det(f))$, using the fixed basis of $V$. Here, $\text{Sym}_n$ is the vector space of symmetric $(n \times n)$-matrices. Note that the image is

$$\mathbb{H}_n \quad := \quad \left\{ (h, z) \in \text{Sym}_n \times \mathbb{C} \,|\, z^2 - \det(h) = 0 \right\}.$$

Since the determinant is an irreducible polynomial, $\mathbb{H}_n$ is a normal variety. The resulting morphism $\text{Hom}(\mathbb{C}^r, V^\vee)/\!/SO(r) \longrightarrow \mathbb{H}_n$ is dominant and birational to $\text{Isom}(\mathbb{C}^r, V^\vee)/SO(r)$. Thus, since both varieties are normal, it is an isomorphism. The above morphism is also equivariant for the remaining $GL(r)$-actions from the left.

If $(\mathcal{A}, \tau)$ is a singular principal $SO(r)$-bundle, the associated section $X \longrightarrow \underline{\text{Hom}}(\mathcal{A}, V^\vee \otimes \mathcal{O}_X)/\!/G$ defines over the open set $U$ where $\mathcal{A}$ is locally free a homomorphism $\iota \colon \mathcal{A}_{|U} \longrightarrow \mathcal{A}^\vee_{|U}$. The associated symmetric pairing $S^2\mathcal{A}_{|U} \longrightarrow \mathcal{O}_U$ can be extended to $\varphi \colon S^2\mathcal{A} \longrightarrow \mathcal{O}_X$, because $\text{codim}(X \setminus U, X) \geq 2$ and $X$ is normal. Moreover, we get an homomorphism $\psi \colon \det \mathcal{A} \longrightarrow \mathcal{O}_X$. Assume $\delta \preceq \delta_\mu$ (cf. Lemma 2.2.3) has degree exactly $\dim X - 1$ and that $(\mathcal{A}, \tau)$ is $\delta$-semistable. Recall that this forces $\mathcal{A}$ to be Mumford-semistable. We claim that $\iota$ must be an isomorphism. Let $j \colon U \longrightarrow X$ be the inclusion. Then, $\iota$ yields

$$\iota_* \colon \mathcal{A} \subset j_*(\mathcal{A}_{|U}) \stackrel{j_*(\iota)}{\longrightarrow} j_*(\mathcal{A}^\vee_{|U}) = \mathcal{A}^\vee.$$

Here, $\mathcal{A} \longrightarrow j_*(\mathcal{A}_{|U})$ is injective, because $\mathcal{A}$ is torsion free, and $\mathcal{A}^\vee \longrightarrow j_*(\mathcal{A}^\vee_{|U})$ is an isomorphism, because $\mathcal{A}^\vee$ is reflexive. Since $\text{codim}(X \setminus U, X) \geq 2$, it follows that $\mathcal{A}^\vee$ is still Mumford-semistable (compare [7], Cor. 3.2.10, p. 67). Both $\det \mathcal{A}$ and $\det \mathcal{A}^\vee$ are trivial. Thus, $\ker \iota_*$ has degree zero. Take the weighted filtration $(\mathcal{A}^\bullet, \underline{\alpha}) := (0 \subset \ker \iota_* \subset \mathcal{A}, (1))$. Then, $M(\mathcal{A}^\bullet, \underline{\alpha})$ has degree at most $\dim X - 2$. On the other hand, one readily checks $\mu(\mathcal{A}^\bullet, \underline{\alpha}; \tau) < 0$, if $\ker \iota_*$ is nontrivial. As $\delta$ has degree exactly $\dim X - 1$, this is impossible. Thus, $(\mathcal{A}, \varphi, \psi)$ is a principal $SO(r)$-sheaf in the sense of [4]. Using [4], Proposition 5.7, one sees that $(\mathcal{A}, \tau)$ will be $\delta$-semistable if and only if $(\mathcal{A}, \varphi, \psi)$ is a semistable principal $SO(r)$-sheaf in the sense of [4]. We infer that we have a set theoretical bijection between the set of isomorphism classes of (semi)stable principal $SO(r)$-sheaves as defined by Gómez and Sols and the set of isomorphism classes of $\delta$-(semi)stable singular $SO(r)$-bundles. Proposition 1.6.1 can be used to identify the moduli functors, using the same argument as above, in order to get an isomorphism between the moduli spaces.

5.3. **The case of the adjoint representation.** Suppose $G'$ is a reductive algebraic group and $\mathfrak{g}$ its Lie algebra. Let $G := G'/\text{Cen}(G')$ be the corresponding semisimple group, and let $\rho \colon G \longrightarrow GL(\mathfrak{g})$ be induced from the adjoint representation of $G$. Then, $G = \text{Aut}(\mathfrak{g})^0$. Here, $\text{Aut}(\mathfrak{g}) \subset GL(\mathfrak{g})$ is the group of those linear transformations respecting the Lie algebra structure of $\mathfrak{g}$. The aim of this section is to show that, for some particular values of the stability parameter $\delta$, our moduli space for $\delta$-semistable singular $G$-bundles will contain only honest singular principal bundles. As Ramanathan's construction shows, it is more convenient to work with $\mathrm{Aut}(\mathfrak{g})$-bundles. In fact, $\mathrm{GL}(\mathfrak{g})/\mathrm{Aut}(\mathfrak{g})$ gets embedded as the $\mathrm{GL}(\mathfrak{g})$-orbit of the Lie bracket on $\mathfrak{g}$ into $\mathrm{Hom}(\mathfrak{g} \otimes \mathfrak{g}, \mathfrak{g})$. Let $\mathcal{V}$ be a vector bundle with fibre $\mathfrak{g}$. Thus, giving a section $X \longrightarrow \underline{\mathrm{Isom}}(\mathfrak{g} \otimes \mathcal{O}_X, \mathcal{V})/\mathrm{Aut}(\mathfrak{g})$ is the same as giving a Lie bracket

$$l \colon \mathcal{V} \otimes \mathcal{V} \longrightarrow \mathcal{V}$$

which induces on every fibre a Lie algebra structure on $\mathfrak{g}$ which is isomorphic to the original one. The map $l$ can also be interpreted as a certain tensor field. This idea can be generalized to map an open set of the parameter space parametrizing honest principal $\mathrm{Aut}(\mathfrak{g})$ bundles to the parameter space of $\delta''$-semisstable tensor fields for a certain polynomial $\delta''$ of degree exactly $\dim X - 1$. If we are able to show that this map is proper, we will be done.

First, let $H$ be a reductive algebraic group acting on a quasi-projective scheme $Y$, and suppose the action is linearized in a line bundle $\mathcal{L}$. Let $G := H^0$ be the connected component of the identity. Then, we have a linearized action of $G$ on $Y$. Recall that the corresponding sets of semistable points are equal for both linearized actions ([11], Prop. 1.15, p. 43).

By what we have just observed, all constructions made so far can also be applied to $\mathrm{Aut}(\mathfrak{g})$. Furthermore, for any singular principal $G$-bundle $(\mathcal{A}, \tau)$, we have an associated singular principal $\mathrm{Aut}(\mathfrak{g})$-bundle $(\mathcal{A}, \tau')$ with

$$\tau' \colon \mathrm{Sym}^*(\mathfrak{g} \otimes \mathcal{A})^{\mathrm{Aut}(\mathfrak{g})} \subset \mathrm{Sym}^*(\mathfrak{g} \otimes \mathcal{A})^G \overset{\tau}{\longrightarrow} \mathcal{O}_X.$$

The above remarks imply indeed that $\tau'$ won't be constant. Moreover, honest singular $G$-bundles are mapped to honest singular $\mathrm{Aut}(\mathfrak{g})$-bundles.

Let $\mathfrak{o} := \#\mathrm{Aut}(\mathfrak{g})/G$. Then, one checks that, for any weighted filtration $(\mathcal{A}^\bullet, \underline{\alpha})$ of $\mathcal{A}$, one has

$$\mu(\mathcal{A}^\bullet, \underline{\alpha}; \tau) \quad = \quad \mathfrak{o} \cdot \mu(\mathcal{A}^\bullet, \underline{\alpha}; \tau').$$

Thus, a $\delta$-(semi)stable singular $G$-bundle leads to an $(\mathfrak{o} \cdot \delta)$-(semi)stable singular $\mathrm{Aut}(\mathfrak{g})$-bundle. The above assignment can clearly be performed also in families. In particular, if $\mathfrak{Y}(G)$ and $\mathfrak{Y}(\mathrm{Aut}(\mathfrak{g}))$ are the respective projective parameter spaces as constructed in the previous chapter, we obtain a finite morphism $\mathfrak{a} \colon \mathfrak{Y}(G)/\!\!/\mathbb{C}^* \longrightarrow \mathfrak{Y}(\mathrm{Aut}(\mathfrak{g}))/\!\!/\mathbb{C}^*$.

Next, we want to construct an open subscheme $\mathfrak{Y}^{\mathrm{h}}(G)$ parametrizing those pairs $([q \colon W \otimes \mathcal{O}_X(-n) \longrightarrow \mathcal{A}], \tau)$ where $H^0(q(n))$ is an isomorphism, $\mathcal{A}$ is torsion free, and $(\mathcal{A}, \tau)$ is an honest singular $G$-bundle. To do so, observe that we can embed $\underline{\mathrm{Hom}}(\mathcal{A}_{\mathfrak{Y}(G)}, \mathfrak{g}^\vee \otimes \mathcal{O}_{\mathfrak{Y}(G) \times X})/\!\!/G$ into a vector bundle $\mathbb{V}_{\mathfrak{Y}(G)}$, and $\mathbb{V}_{\mathfrak{Y}(G)}$ into the projective bundle $\mathbb{P}(\mathbb{V}^\vee_{\mathfrak{Y}(G)} \oplus \mathcal{O}_{\mathfrak{Y}(G) \times X})$. Let $U_{\mathrm{lf}}$ be the maximal open subset of $\mathfrak{Y}(G) \times X$ where $\mathcal{A}_{\mathfrak{Y}}$ is locally free. Then, we have the open subscheme

$$\underline{\mathrm{Isom}}(\mathcal{A}_{\mathfrak{Y}(G)|U_{\mathrm{lf}}}, \mathfrak{g}^\vee \otimes \mathcal{O}_{U_{\mathrm{lf}}}) \subset \mathbb{P}(\mathbb{V}^\vee_{\mathfrak{Y}(G)|U_{\mathrm{lf}}} \oplus \mathcal{O}_{U_{\mathrm{lf}}}).$$

Let $\mathbb{I}_{U_{\mathrm{lf}}}$ be its complement and $\mathbb{I}_{\mathfrak{Y}(G)} \subset \mathbb{P}(\mathbb{V}^\vee_{\mathfrak{Y}(G)} \oplus \mathcal{O}_{\mathfrak{Y}(G) \times X})$ the closure of $\mathbb{I}_{U_{\mathrm{lf}}}$. By elimination, the image $\mathfrak{I}_{\mathfrak{Y}(G)}$ of $\mathbb{I}_{\mathfrak{Y}(G)}$ is a closed subscheme of $\mathfrak{Y}(G)$. Then,

$$\mathfrak{Y}^{\mathrm{h}}(G) \quad := \quad \mathfrak{Y}^0(G) \setminus \mathfrak{I}_{\mathfrak{Y}(G)}$$

is the parameter space we are looking for. Similarly, we define $\mathfrak{Y}^{\mathrm{h}}(\mathrm{Aut}(\mathfrak{g}))$. Note that $\mathfrak{Y}^{\mathrm{h}}(G)/\!/\mathbb{C}^* = \mathfrak{Y}(G)/\!/\mathbb{C}^* \times_{\mathfrak{a}} \left( \mathfrak{Y}^{\mathrm{h}}(\mathrm{Aut}(\mathfrak{g}))/\!/\mathbb{C}^* \right)$, so that we have a finite morphism

$$\mathfrak{a}^{\mathrm{h}} \colon \mathfrak{Y}^{\mathrm{h}}(G) \longrightarrow \mathfrak{Y}^{\mathrm{h}}(\mathrm{Aut}(\mathfrak{g})).$$

Next, we can embed the $\mathrm{GL}(\mathfrak{g})$-module $\mathrm{Hom}(\mathfrak{g}, \mathfrak{g} \otimes \mathfrak{g})$ into

$$V_{a,b,c} := \big(\mathfrak{g}^{\otimes a} \otimes (\overset{\dim \mathfrak{g}}{\bigwedge} \mathfrak{g})^{\otimes -b}\big)^{\oplus c}$$

for some non-negative integers $a, b, c$. Let $\mathbb{S} \subset V_{a,b,c}$ be the submodule supplementary to $\mathrm{Hom}(\mathfrak{g}, \mathfrak{g} \otimes \mathfrak{g})$. We use the Killing form of $\mathfrak{g}$ to identify the $\mathrm{GL}(\mathfrak{g})$-modules $\mathrm{Hom}(\mathfrak{g}, \mathfrak{g}^\vee)$ and $\mathrm{Hom}(\mathfrak{g}, \mathfrak{g})$. The immersion $\mathrm{GL}(\mathfrak{g})/\mathrm{Aut}(\mathfrak{g}) \longrightarrow \mathrm{Hom}(\mathfrak{g} \otimes \mathfrak{g}, \mathfrak{g})$ provides a $\mathrm{GL}(\mathfrak{g})$-equivariant rational map

$$
\begin{aligned}
r \colon \mathbb{P}(\mathrm{Hom}(\mathfrak{g}, \mathfrak{g}^\vee)^\vee)/\!/\mathrm{Aut}(\mathfrak{g}) \quad &\cong \quad \mathbb{P}(\mathrm{Hom}(\mathfrak{g}, \mathfrak{g}))/\!/\mathrm{Aut}(\mathfrak{g}) \\
&\dashrightarrow \quad \mathbb{P}(\mathrm{Hom}(\mathfrak{g}, \mathfrak{g} \otimes \mathfrak{g})) \hookrightarrow \mathbb{P}(V_{a,b,c}).
\end{aligned}
$$

Let $\mathcal{N}$ be the pullback of $\mathcal{O}(1)$ to $\mathbb{P}(\mathrm{Hom}(\mathfrak{g}, \mathfrak{g}^\vee)^\vee)/\!/\mathrm{Aut}(\mathfrak{g})$ via this rational map and let $\mathcal{M}_s$ be the line bundle on $\mathbb{P}(\mathrm{Hom}(\mathfrak{g}, \mathfrak{g}^\vee)^\vee)/\!/\mathrm{Aut}(\mathfrak{g})$ coming from its embedding into $\mathbb{P}(\mathbb{U}(s))$ by pullback of $\mathcal{O}(1)$. Then, for some natural numbers $a$ and $b$, we have $\mathcal{M}_s^{\otimes a} \cong \mathcal{N}^{\otimes b}$ as $\mathrm{SL}(\mathfrak{g})$-linearized line bundles. To see this, it is enough to prove the analogous assertion for the $(\mathrm{SL}(\mathfrak{g}) \times \mathrm{Aut}(\mathfrak{g}))$-linearized line bundles on $\mathbb{P}(\mathrm{Hom}(\mathfrak{g}, \mathfrak{g}^\vee)^\vee)$ obtained by pullback via the quotient map. For these, the claim follows from the fact that any two $(\mathrm{SL}(\mathfrak{g}) \times \mathrm{Aut}(\mathfrak{g}))$-linearizations in a given line bundle differ by a character and every character of $\mathrm{SL}(\mathfrak{g}) \times \mathrm{Aut}(\mathfrak{g})$ is obviously torsion. Thus, for any one parameter subgroup $\lambda$ of $\mathrm{SL}(\mathfrak{g})$ and any point $[f] \in \mathbb{P}(\mathrm{Hom}(\mathfrak{g}, \mathfrak{g}^\vee)^\vee)$ in which $r$ is defined, we have

$$(1) \qquad\qquad as! \cdot \mu(\lambda, [f]) \quad \geq \quad b \cdot \mu(\lambda, r([f])).$$

For fixed $s$, the ratio $a/b$ is well defined, and for $\mathcal{M}_{s+1}$ we have $\mathcal{M}_{s+1}^{\otimes a/(s+1)} = \mathcal{M}_s^{\otimes a} \cong \mathcal{N}^{\otimes b}$, so that the same goes for the ratio $\eta := as!/b$.

Suppose $(\mathcal{A}, \tau')$ is an honest singular $\mathrm{Aut}(\mathfrak{g})$-bundle. Let $U$ be the maximal open set over which $\mathcal{A}$ is locally free. The restriction of $\tau'$ to $U$ may be interpreted as a Lie bracket

$$\widetilde{l}_U \colon \mathcal{A}_U^\vee \otimes \mathcal{A}_U^\vee \longrightarrow \mathcal{A}_U^\vee.$$

We use Ramanathan's morphism

$$\mathrm{Ra} \colon \mathrm{Hom}(\mathcal{A}_U^\vee, \underline{\mathrm{End}}(\mathcal{A}_U^\vee)) \quad \longrightarrow \quad \mathrm{Hom}(S^2\mathcal{A}_U^\vee, \mathcal{O}_U)$$

with

$$\mathrm{Ra}(l) \colon \mathcal{A}_U^\vee \otimes \mathcal{A}_U^\vee \overset{l \otimes l}{\longrightarrow} \underline{\mathrm{End}}(\mathcal{A}_U^\vee) \otimes \underline{\mathrm{End}}(\mathcal{A}_U^\vee) \overset{\mathrm{mult}}{\longrightarrow} \mathrm{End}(\mathcal{A}_U^\vee) \overset{\mathrm{trace}}{\longrightarrow} \mathcal{O}_U.$$

The bilinear form $\mathrm{Ra}(\widetilde{l}_U)$ is non-degenerate, because it is fibrewise the Killing form. The induced isomorphism $\mathcal{A}_U \longrightarrow \mathcal{A}_U^\vee$ provides $\mathcal{A}_U$ itself with a Lie bracket

$$l'_U \colon \mathcal{A}_U \otimes \mathcal{A}_U \longrightarrow \mathcal{A}_U$$

or, equivalently, a trilinear map

$$l_U \colon \mathcal{A}_U \otimes \mathcal{A}_U \otimes \mathcal{A}_U^\vee \longrightarrow \mathcal{O}_U.$$

Since $\det \mathcal{A}_U$ is trivial, the projection of $V_{a,b,c}$ onto $\mathrm{Hom}(\mathfrak{g}, \mathfrak{g} \otimes \mathfrak{g})$ together with $l_U$ induces a tensor field

$$t_U \colon \mathcal{A}_U^{\otimes a \oplus c} \longrightarrow \mathcal{O}_U$$

which can be extended to

$$t\colon \mathcal{A}^{\otimes a\oplus c} \longrightarrow \mathcal{O}_X.$$

Let $\delta''$ be a positive polynomial of degree at most $\dim X - 1$. If $(\mathcal{A}, t)$ is $\delta''$-(semi)stable tensor field, then $(\mathcal{A}, \tau')$ is $(\eta \cdot \delta'')$-(semi)stable, by equation (1). Observe that the analogon of Lemma 2.2.3 ii) holds also for tensor fields. We now assume that $\delta''$ has degree exactly $\dim X - 1$ and that for every $\delta''$-semistable tensor field $(\mathcal{A}, t)$ the sheaf $\mathcal{A}$ is Mumford-semistable.

Next, let $\mathfrak{S}$ be the projective parameter space of equivalence classes of pairs $([q\colon W \otimes \mathcal{O}_X(-n) \longrightarrow \mathcal{A}], [t])$ where $q\colon W \otimes \mathcal{O}_X(-n) \longrightarrow \mathcal{A}$ is a quotient with Hilbert polynomial $P$ and $t\colon \mathcal{A}^{\otimes a\oplus c} \longrightarrow \mathcal{O}_X$ is non-trivial tensor field, and the equivalence relation is the obvious one.

The passage from $\tau'$ to the tensor field $t$ can be clearly performed in families, too, and, thus, there is a morphism

$$\mathfrak{s}\colon \mathfrak{Y}^h(G)/\!\!/\mathbb{C}^* \longrightarrow \mathfrak{Y}^h(\mathrm{Aut}(\mathfrak{g}))/\!\!/\mathbb{C}^* \longrightarrow \mathfrak{S}.$$

Let $\mathfrak{S}^{\delta''-ss}$ be the open subscheme, parametrizing $([q\colon W \otimes \mathcal{O}_X(-n) \longrightarrow \mathcal{A}], [t])$ where $H^0(q(n))$ is an isomorphism, $\mathcal{A}$ is torsion free, and $(\mathcal{A}, t)$ is a $\delta''$-semistable tensor field. Our efforts finally culminate in

**Theorem 5.3.1.** *The induced map*

$$\mathfrak{s}^0\colon \mathfrak{s}^{-1}\big(\mathfrak{S}^{\delta''-ss}\big) \longrightarrow \mathfrak{S}^{\delta''-ss}$$

*is proper.*

This theorem implies that the set $\mathfrak{s}^{-1}\big(\mathfrak{S}^{\delta''-ss}\big)$ is a saturated open subset of the $\mathrm{SL}(W)$-semistable points in $\mathfrak{Y}(G)/\!\!/\mathbb{C}^*$. Thus,

$$\mathcal{HSPB}_P^{\delta-ss} \quad := \quad \mathfrak{s}^{-1}\big(\mathfrak{S}^{\delta''-ss}\big)/\!\!/\mathrm{SL}(W)$$

is an open subscheme of $\mathfrak{M}(\rho)_P^{\delta-ss}$, $\delta := (\eta/\mathfrak{o}) \cdot \delta''$. On the other hand, $\mathcal{HSPB}_P^{\delta-ss}$ maps properly to the moduli scheme of $\delta''$-semistable tensor fields. The latter space being projective, we conclude that $\mathcal{HSPB}_P^{\delta-ss}$ is projective and a union of connected components of $\mathfrak{M}(\rho)_P^{\delta-ss}$. The space $\mathcal{HSPB}_P^{\delta-ss}$ supplies a natural compactification of Hyeon's moduli space which is analogous to the Gieseker compactification of the moduli space of Mumford-stable vector bundles.

*Proof of Theorem 5.3.1.* We will apply the valuative criterion for properness. Observe that, by the properness of the map $\mathfrak{o}$, it suffices to prove the analogous assertion for the space $\mathfrak{Y}^h(\mathrm{Aut}(\mathfrak{g}))$. Let $(C, 0)$ be the spectrum of DVR $R$ with quotient field $K$. Suppose we have a morphism $h\colon C \longrightarrow \mathfrak{S}^{\delta''-ss}$ which lifts over $\mathrm{Spec}\,K$ to $\mathfrak{Y}^h(\mathrm{Aut}(\mathfrak{g}))/\!\!/\mathbb{C}^*$. By Luna's étale slice theorem, we may, after a finite extension of $K$, assume that it even lifts to $\mathfrak{Y}^h(\mathrm{Aut}(\mathfrak{g}))$.

By pullback of the universal family, we obtain a family $(q_C\colon W \otimes \pi_X^* \mathcal{O}_X(-n) \longrightarrow \mathcal{A}_C, t_C\colon \mathcal{A}^{\otimes a\oplus c} \longrightarrow \mathcal{O}_{C\times X})$. Over $C \setminus \{0\}$, the map $h$ lifts to the scheme $\mathfrak{Y}^h(\mathrm{Aut}(\mathfrak{g}))$. This lifting is induced by a family $(q_K'\colon W \otimes \pi_X^* \mathcal{O}_X(-n) \longrightarrow \mathcal{A}_K', \tau_K')$. This family defines a family $(q_K', t_K')$ of $\delta''$-semistable tensor fields which is equivalent to the family $(q_{C|\mathrm{Spec}\,K\times X}, t_{C|\mathrm{Spec}\,K\times X})$ in the sense that there are an isomorphism $\psi_K\colon \mathcal{A}_K' \longrightarrow \mathcal{A}_{C|\mathrm{Spec}\,K\times X}$ and $\lambda \in K^*$ with $q_{C|\mathrm{Spec}\,K\times X} = \psi_K \circ q_K'$ and $\lambda \cdot t_K' = t_{C|\mathrm{Spec}\,K\times X} \circ \psi_K^{\otimes a\oplus c}$. We can read this backwards, too. Let $\mathfrak{U}_{1f} \subset \mathrm{Spec}\,K \times X$ be the open set where $\mathcal{A}_K'$ and $\mathcal{A}_{C|\mathrm{Spec}\,K\times X}$ are locally free. The tensor fields $t_{K|\mathfrak{U}_{1f}}'$ and $t_{C|\mathfrak{U}_{1f}}$ factorize over $\mathcal{A}_{K|\mathfrak{U}_{1f}}' \otimes \mathcal{A}_{K|\mathfrak{U}_{1f}}' \otimes \mathcal{A}_{K|\mathfrak{U}_{1f}}'^\vee$ and $\mathcal{A}_{C|\mathfrak{U}_{1f}} \otimes \mathcal{A}_{C|\mathfrak{U}_{1f}} \otimes \mathcal{A}_{C|\mathfrak{U}_{1f}}^\vee$, respectively. Let $l_K'$ and

$l_K$ be the associated Lie brackets. We clearly have $\lambda^2 \cdot \mathrm{Ra}(l_K') = \mathrm{Ra}(l_K) \circ S^2 \psi_K$, in particular, $\mathrm{Ra}(l_K)$ is non-degenerate. The induced isomorphism $\mathcal{A}_{C|\mathcal{U}_{\mathrm{lf}}} \longrightarrow \mathcal{A}_{C|\mathcal{U}_{\mathrm{lf}}}^\vee$ can be used to equip $\mathcal{A}_{C|\mathcal{U}_{\mathrm{lf}}}^\vee$ with a Lie bracket which comes from a section $\mathcal{U}_{\mathrm{lf}} \longrightarrow \underline{\mathrm{Isom}}(\mathcal{A}_{C|\mathcal{U}_{\mathrm{lf}}}, \mathfrak{g}^\vee \otimes \mathcal{O}_{\mathcal{U}_{\mathrm{lf}}})$. By Proposition 1.6.1, this extends to a family $(q_{C|\mathrm{Spec}\,K \times X}, \tau_K'')$ which defines a map to $\mathfrak{Y}^h(\mathrm{Aut}(\mathfrak{g}))/\!/\mathbb{C}^*$ which is, of course, the lifting of $h_{|\mathrm{Spec}\,K}$ we started with.

Now, let $\iota\colon \mathcal{U}_{\mathrm{lf}} \hookrightarrow C \times X$ stand for the maximal open subset where $\mathcal{A}_C$ is locally free. Let $\mathcal{A}_{\mathbb{S}, \mathcal{U}_{\mathrm{lf}}}$ be the vector bundle with fibre $\mathbb{S}$ associated with $\mathcal{A}_{\mathcal{U}_{\mathrm{lf}}}$. Then, $\iota_* t_{C|\mathcal{U}_{\mathrm{lf}}}$ vanishes on $\iota_* \mathcal{A}_{\mathbb{S}, \mathcal{U}_{\mathrm{lf}}}$, because this is a closed condition, by Proposition 1.5.1, which holds on $\mathcal{U}_{\mathrm{lf}} \cap (\mathrm{Spec}\,K \times X)$. In other words, on $\mathcal{U}_{\mathrm{lf}}$ we have a bilinear map

$$l_C \colon \mathcal{A}_{\mathcal{U}_{\mathrm{lf}}} \otimes \mathcal{A}_{\mathcal{U}_{\mathrm{lf}}} \longrightarrow \mathcal{A}_{\mathcal{U}_{\mathrm{lf}}}$$

which satisfies the Jacobi identity, by continuity. Now, we continue as Ramanathan [13], page 445. First, we claim that $\mathrm{Ra}(l_{C|U})$, $U := \mathcal{U}_{\mathrm{lf}} \cap \{0\} \times X$, is non-degenerate. If this were not the case, there would be a non-trivial saturated subsheaf $\mathcal{F}$ of $\mathcal{A}_{C|\{0\} \times X}$ of degree zero such that $\mathcal{F}_{|U}$ lies in the center of $l_{C|U}$. If we take the weighted filtration $(0 \subset \mathcal{F} \subset \mathcal{A}, (1))$, we easily compute $\mu((0 \subset \mathcal{F} \subset \mathcal{A}, (1)); t_{C|\{0\} \times X}) = \mathrm{rk}\,\mathcal{F} - \mathrm{rk}\,\mathcal{A}$. As $M(0 \subset \mathcal{F} \subset \mathcal{A}, (1))$ has degree at most $\dim X - 2$, because $\deg \mathcal{F} = \deg \mathcal{A} = 0$, the polynomial

$$M(0 \subset \mathcal{F} \subset \mathcal{A}, (1)) + \delta'' \cdot \mu\big((0 \subset \mathcal{F} \subset \mathcal{A}, (1)); t_{C|\{0\} \times X}\big)$$

would be stricly negative, a contradiction. Thus, $l_C$ induces on every fibre the structure of a semisimple Lie algebra which, by rigidity, is isomorphic to $\mathfrak{g}$. By the construction we have described before and equation (1), we get a family $(q_C, \tau_C')$ of $(\eta \cdot \delta'')$-semistable honest $\mathrm{Aut}(\mathfrak{g})$-bundles. This family defines a morphism $k\colon C \longrightarrow \mathfrak{Y}^h(\mathrm{Aut}(\mathfrak{g}))$, and we have seen before that $k_{|\mathrm{Spec}\,K}$ is just the lifting of $h$ we started with, so that we are done. $\qquad\square$

*The work of Gómez and Sols.* In [5], which appeared after the first version of my paper had been finished, Gómez and Sols have announced the construction of moduli spaces for "principal $G'$-sheaves". These are objects $(\mathcal{E}, \varphi, \xi)$, consisting of a torsion free sheaf $\mathcal{E}$ with fibre $\mathfrak{g}^c$, where $\mathfrak{g}^c$ is the Lie algebra of $[G, G]$, a bilinear pairing $\varphi\colon \mathcal{E} \otimes \mathcal{E} \longrightarrow \mathcal{E}^{\vee\vee}$ which induces on every fibre over the maximal open subset $U$ where $\mathcal{E}$ is locally free a Lie algebra structure isomorphic to $\mathfrak{g}^c$, and a lifting $\xi$ of the $\mathrm{Aut}(\mathfrak{g}^c)$-bundle on $U$ defined by $(\mathcal{E}, \varphi)$ to a $G'$-bundle. The datum $(\mathcal{E}, \varphi)$ is obviously equivalent to the datum of an honest $\mathrm{Aut}(\mathfrak{g})$-bundle.

For pairs $(\mathcal{E}, \varphi)$, where $\varphi$ can now be any bilinear pairing, one has again a natural notion of $\delta$-semistability. For $(\mathcal{E}, \varphi)$ coming from an honest $\mathrm{Aut}(\mathfrak{g})$-bundle, this notion of $\delta$-semistability is the same as the $\delta$-semistability of the associated tensor field $(\mathcal{E}, t)$ which we used in our proof.

Gómez and Sols call $(\mathcal{E}, \varphi, \xi)$ *semistable*, if $\mathcal{E}$ is a Gieseker-semistable sheaf. If $\delta$ is just a small constant, then $\delta$-semisatbility of $(\mathcal{E}, \varphi)$ implies that $\mathcal{E}$ is itself a Gieseker semistable sheaf. In [4], it is asserted first that, for all polynomial $\delta$, there is a projective moduli space for $\delta$-semistable pairs $(\mathcal{E}, \varphi)$ and that this can be used to construct a projective moduli space for principal $G'$-sheaves. If we take $G' = G$, then this statement is clearly equivalent to the statement that $\mathfrak{s}$ is also proper for $\delta$ a small constant. I have not yet been able to verify this.

Universität GH Essen, FB6 Mathematik & Informatik, D-45117 Essen, Germany
*E-mail address:* `alexander.schmitt@uni-essen.de`